\documentstyle[12pt,righttag]{amsart}

\numberwithin{equation}{section}
\newtheorem{thm}{Theorem}[section]

\newtheorem{lem}[thm]{Lemma}
\newtheorem{prop}[thm]{Proposition}
\newtheorem{conj}[thm]{Conjecture}

\theoremstyle{definition}
\newtheorem{defn}{Definition}[section]
\theoremstyle{remark}
\newtheorem{rem}{Remark}
\newtheorem{exm}[rem]{Example}
\renewcommand{\therem}

\title[Skew Young diagram method]
{Skew Young diagram method in spectral decomposition of
integrable lattice models}

\begin{document}
\newcommand{\ignoreit}[2]{\null}

\maketitle
\begin{center}
{\sc Anatol N.\ Kirillov\footnote[1]{
Permanent address: Steklov Mathematical Institute,
Fontanka 27, St.\ Petersburg, 191011, Russia},
Atsuo Kuniba$^2$, and Tomoki Nakanishi\footnote[3]{
Permanent address: Department of Mathematics,
Nagoya University, Chikusa-ku, Nagoya 464, Japan}}\\
 ~\\
 ~\\
{\it 
$^1$Department of Mathematical Sciences, 
University of Tokyo\\
Komaba, Meguro-ku, Tokyo 153, Japan\\
 ~\\
$^2$Institute of Physics,
University of Tokyo\\
Komaba, Meguro-ku, Tokyo 153, Japan\\
 ~\\
$^3$Department of Mathematics,
University of North Carolina\\
Chapel Hill, NC 27599, USA}\\
~\\
July, 1996
\end{center}

\begin{abstract}
The spectral decomposition of the path space of the vertex model associated
to the vector representation of the quantized affine algebra
$U_q(\widehat{sl}_n)$ is studied. We give a one-to-one correspondence between
the spin configurations and the semi-standard tableaux of skew Young diagrams.
As a result we obtain a formula of the characters for the degeneracy of the
spectrum  in terms of skew Schur functions. We conjecture that our
result describes the $sl_n$-module contents
of the Yangian $Y(sl_n)$-module structures of the
level 1 integrable modules of the affine Lie algebra $\widehat{sl}_n$.
An analogous result is obtained also for a vertex model
associated to the quantized twisted affine
algebra $U_q(A^{(2)}_{2n})$, where $Y(B_n)$ characters
appear for the degeneracy of the spectrum.
The relation to the spectrum of the Haldane-Shastry and 
the Polychronakos
models are also discussed.
\end{abstract}

\newpage

\section{Introduction}

The corner transfer matrix (CTM) has been attracting much
attention in the recent study of integrable lattice models
based on the Yang-Baxter equation \cite{baxter}.
The CTM acts on the space of paths, which
is often identified with the semi-infinite
tensor product of a finite-dimensional
quantum group module.
It is well known that the trace of the CTM
of a vertex model associated to the
quantized affine algebra $U_q(\widehat{g})$ is related to
affine Lie algebra characters.
We call this correspondence the DJKMO 
(Date-Jimbo-Kuniba-Miwa-Okado) correspondence.

In \cite{arakawa} a fine structure of the CTM
spectrum, called the spectral decomposition, is studied
in the $U_q(\widehat{sl}_2)$ vertex models.
The idea behind it is as follows:
The logarithmic derivative of the CTM
is regarded as the energy operator, or the
Hamiltonian, of the
path space.
As a nature of an integrable system,
we expect that there exists a family of
commuting operators
(the integrals of motion)
which act on the path space and commute
with the Hamiltonian.
The spectral decomposition is  the simultaneous
diagonalization of these integrals of motion.
The degeneracy of the spectrum, then, reflects the
non-abelian symmetry which commutes with these integrals of motion.

In the meanwhile, the action of the Yangian algebra $Y(sl_2)$
is defined on the level $1$ integrable modules of the untwisted
affine Lie algebra $\widehat{sl}_2$ in \cite{haldanePRL},
and their $Y(sl_2)$-module structures are determined \cite{bernardNP}.
It turns out that
the degeneracy of the
spectrum of the CTM Hamiltonian in the $U_q(\widehat{sl}_2)$ vertex
 model exactly describes the $Y(sl_2)$-module structure
of the level $1$ integrable modules \cite{arakawa}.

We believe this coincidence is a universal phenomena.
Namely, we expect that
a similar coincidence occurs between the
spectrum of the commutant of
the non-abelian, probably
some quantum group, symmetry of a conformal field theory
and the  CTM spectrum  in the corresponding lattice model.

Motivated by this expectation, in this paper we study the spectral
decomposition of the vertex model of the vector
representation of the $U_q(\widehat{sl}_n)$. The counter part of the
DJKMO correspondence is the level $1$ integrable modules of 
$\widehat{sl}_n$.
The action of the  Yangian $Y(sl_n)$ on these
modules is defined in \cite{schoutens}, but
the $Y(sl_n)$-module structure is not fully studied yet
(however, see \cite{bernardJP} for a related result).
In this paper we determine the characters of
the degeneracy of the spectrum, and show that they are the characters
of irreducible $Y(sl_n)$-modules  as expected.
Therefore, we conjecture that our spectral decomposition
exactly describes the $Y(sl_n)$-module structure of the
level $1$ integrable modules at the character level.

Conceptually, the spectral decomposition in the $sl_n$ case
is formulated just as in the case of $sl_2$.
However, due to the complexity of the irreducible
modules of $Y(sl_n)$ for $n\geq 3$, the incidence
matrix technique used in \cite{arakawa} is not
very efficient.
A key  to overcome this difficulty is the observation
that
there exists a natural one-to-one correspondence
between the paths of the vertex model
and the semi-standard tableaux of certain skew Young
diagrams. 
The appearance of the skew Young diagrams is not
quite unexpected, because they label a
family of $Y(sl_n)$-modules, called the tame modules
in \cite{nazarov}.
Thanks to this correspondence, we can express the
characters of the degeneracy of the spectrum
in terms of the skew Schur functions. This enables
us to identify them as irreducible $Y(sl_n)$
characters.

It is possible to extend our analysis for
a vertex model associated to the quantized
twisted affine algebra $U_q(A^{(2)}_{2n})$.
The characters of the degeneracy of the spectrum
are analogues of the skew Schur functions.
They are conjectured in \cite{kunibas}
to be irreducible $Y(B_n)$ characters.

The content of the paper is as follows:
In section 2 we review the DJKMO correspondence
for the vertex model in the $sl_n$ case. In section 3 we formulate
the spectral decomposition of the model.
In section 4 some properties of skew Young diagrams and 
skew Schur functions are given.
In section 5 we describe the correspondence between
the configurations of the vertex model and the
semi-standard tableaux, and determine the
characters of the degeneracy of the spectrum.
In section 6 the identifications with
irreducible $Y(gl_n)$ and $Y(sl_n)$ characters
are given.
In section 7 we discuss the relation with
the spectrum of other type of spin models, such as
the Haldane-Shastry model and the Polychronakos
model. 
In section 8 an analogous result for an $U_q(
A^{(2)}_{2n})$ vertex model is presented.
In Appendix A the equality between the skew
Schur functions and the Rogers-Szeg\"o polynomials
is proved. In Appendix B a new combinatorial description
of the Kostka-Foulkes polynomials is given.
In Appendix C we describe the level $1$ character of $A^{(2)}_{2n}$.

\section{DJKMO correspondence}

We review the correspondence between the CTM spectrum of the
vertex models of the vector representation of $U_q(\widehat{sl}_n)$
and the affine Lie algebra characters of $\widehat{sl}_n$ \cite{DJKMO}.

For given two infinite sequences,
 $\vec{a}=(a_1,a_2,\dots)$
and 
$\vec{b}=(b_1,b_2,\dots)$, of any kind of
objects $a_i,b_i$,
we write $\vec{a}\approx \vec{b}$
if $a_i\neq b_i$ only
 for finitely many $i$.
We often use a shorthand notation
$\vec{a}=(a_1,\dots,a_k,(a_{k+1},\dots,a_{k+m})^\infty)$
for such a periodic sequence as
$\vec{a}=(a_1,\dots,a_k,a_{k+1},
\dots,a_{k+m},a_{k+1},\dots,a_{k+m},\dots)$.

\par
Let $\overline{\Lambda}_1,\dots,\overline\Lambda_{n-1}$
 be the fundamental weights
of the Lie algebra $sl_n$, and let
\begin{equation}\label{eq:weights}
\epsilon_i=\overline\Lambda_i-\overline\Lambda_{i-1}
\end{equation}
for $i=1,\dots,n$ 
with $\overline\Lambda_0=\overline\Lambda_n=0$.
Then  $B(\overline\Lambda_1)
=\{\epsilon_1,\dots,\epsilon_n\}$ is the set of
all the weights of
the irreducible
representation (vector representation)
of  $sl_n$ whose highest weight is $\overline\Lambda_1$.
 We give a total ordering in $B(\overline\Lambda_1)$ as
$\epsilon_1 \prec \epsilon_2 \prec \dots \prec \epsilon_{n}$.

\par
We define the {\it local energy  function}
$H: B(\overline\Lambda_1)
\times B(\overline\Lambda_1)\rightarrow \{0,1\}$ as
\begin{equation}\label{eq:localenergy}
H(\epsilon_i,\epsilon_j)=
\begin{cases}
0 & \mbox{if $\epsilon_i \prec \epsilon_j$}, \cr
1 & \mbox{if $\epsilon_i\succeq \epsilon_j$}. \cr
\end{cases}
\end{equation}
The function $H$ will play an essential role in our study.
It is  the logarithm of the $R$-matrix
associated to the vector representation of $U_q(\widehat{sl}_n)$
in the limit $q\rightarrow 0$.

\par
An infinite sequence $\vec{s}=(s_i)$, 
 $s_i\in B(\overline\Lambda_1)$,
 is called a {\it spin configuration} if
it has a form
\begin{equation*}
\vec{s}=(s_1,\dots,s_m,(\epsilon_1,\epsilon_2,\dots,\epsilon_n)^\infty),
\end{equation*}
where $s_1,\dots,s_m$ is an arbitrary finite sequence.
Equivalently, $\vec{s}$ is a spin configuration if
$\vec{s}\approx \vec{s}^{(k)}$ for some $k=0,\dots,n-1$,
where
\begin{equation}\label{eq:boundary}
\vec{s}^{(k)}
=(\epsilon_{1},\epsilon_{2},
\dots,\epsilon_k,(\epsilon_1,\epsilon_2,\dots,\epsilon_n)^\infty).
\end{equation}
The set of all the spin configurations 
${\cal S}$ has a natural decomposition
\begin{equation}
{\cal S}=\bigsqcup_{k=0}^{n-1} {\cal S}^{(k)},
\quad
{\cal S}^{(k)}=\{ \vec{s}\mid \vec{s}\approx \vec{s}^{(k)}\}.
\end{equation}

\par
For $\vec{s}=(s_i)\in {\cal S}^{(k)}$ we define
its {\it energy} $E(\vec{s})$ and $sl_n$-{\it weight}
${\rm wt}(\vec{s})$ as
\begin{align}\label{eq:energyweight}
E(\vec{s})&=\sum_{i=1}^\infty i \left\{H(s_i,s_{i+1})-H(s^{(k)}_i,
s^{(k)}_{i+1})\right\},
\tag{\theequation a}
\\
{\rm wt}(\vec{s})&=\overline\Lambda_k
+\sum_{i=1}^\infty \left(s_i - s_i^{(k)}
\right).
\tag{\theequation b}
\end{align}
\addtocounter{equation}{1}

\begin{prop}\label{prop:weight}
Let $\vec{s}=(s_1,\dots,s_{m},
(\epsilon_1,\dots,\epsilon_n)^\infty)$
be any element of ${\cal S}^{(k)}$.
Then
\begin{equation*}
{\rm wt}(\vec{s})=\sum_{i=1}^{m} s_i.
\end{equation*}
\end{prop}
\begin{pf}
Since $m\equiv k$ modulo $n$,
$\sum_{i=1}^m s^{(k)}_i=\overline\Lambda_k$.
Thus
\begin{equation*}
\sum_{i=1}^{m} s_i =
\overline\Lambda_k + \sum_{i=1}^{m}
(s_i-s_i^{(k)})
=
\overline\Lambda_k + \sum_{i=1}^{\infty}
(s_i-s_i^{(k)})
={\rm wt}(\vec{s}).
\end{equation*}
\end{pf}

There is a remarkable connection between
the partition function of ${\cal S}^{(k)}$
and an affine Lie algebra character.

\begin{thm}[DJKMO correspondence \cite{DJKMO,kang}] \label{thm:DJKMO}
For $k=0,1,\dots,n-1$, let ${\cal L}(\Lambda_k)$
be the level 1 integrable module of
the untwisted affine Lie algebra  $\widehat{sl}_n$ whose 
highest weight is the $k${\rm th} fundamental 
weight $\Lambda_k$ of $\widehat{sl}_n$.
Then the following equality holds:
\begin{align}\label{eq:djkmo}
{\rm ch}\, {\cal L}(\Lambda_k)
&=
q^{\Delta_k - c/24}
\sum_{\vec{s}\in {\cal S}^{(k)}}
 q^{E(\vec{s})}
e^{{\rm wt}(\vec{s})}
\tag{\theequation a}\\
&=q^{\Delta_k - c/24}
\sum_{\vec{s}\in {\cal S}^{(n-k)}}
 q^{E(\vec{s})}
e^{-{\rm wt}(\vec{s})},
\tag{\theequation b}
\end{align}
\noindent
where
\addtocounter{equation}{1}
${\rm ch}\, {\cal L}(\Lambda_k)$,
$\Delta_k=k(n-k)/2n$, and $c=n-1$ are
the (normalized) character,
the conformal dimension, and the Virasoro central charge
of ${\cal L}(\Lambda_k)$, respectively \cite{kac}.
In (\ref{eq:djkmo}b) ${\cal S}^{(n)}
={\cal S}^{(0)}$.
\end{thm}

\begin{rem}
Often spin configurations are described by an alternative notion,
{\it paths}.
An infinite sequence $\vec{p}=(p_i)$ of
$\widehat{sl}_n$-weights is called a path if it satisfies the
conditions:
(i) $p_{i+1}-p_i \in B(\overline\Lambda_1)$,
(ii) $\vec{p}=(p_1,\dots,p_m,
(\Lambda_0,\Lambda_1,\dots,
\Lambda_{n-1})^\infty)$.
It is clear that the map
$\vec{p}=(p_i)\mapsto
\vec{s}=(p_{i+1}-p_i)$
is a bijection
from the set of all the
 paths ${\cal P}$ to ${\cal S}$.
By  this identification our spin configurations are also
called paths in some literature.
The $sl_n$-weight of a path $\vec{p}=(p_i)$
is defined as ${\rm wt}(\vec{p})=\overline{p_1}$.
Then, under the bijection $\vec{p}\mapsto \vec{s}$,
${\rm wt}(\vec{p})=-{\rm wt}(\vec{s})$ holds.
It is this context where (\ref{eq:djkmo}b)
is proved in \cite{DJKMO}.
The expression (\ref{eq:djkmo}a) follows from (\ref{eq:djkmo}b)
by the Dynkin diagram automorphism $\alpha_i
\leftrightarrow \alpha_{n-i}$ of $sl_n$.
\end{rem}

\section{Spectral decomposition}

We introduce the {\it local energy map} $h:{\cal S}\rightarrow
\{ 0,1\}^{\bf N}$ such that
\begin{equation}\label{eq:maph}
h:\vec{s}=(s_i)\mapsto \vec{h}=(h_i),\quad
h_i=H(s_i,s_{i+1}).
\end{equation}
Each number $h_i$ is called the {\it $i${\rm th}
 local energy} of $\vec{s}$.
We call the image ${\rm Sp}=h({\cal S})$
the {\it spectrum} of ${\cal S}$. Let
${\rm Sp}^{(k)}=h({\cal S}^{(k)})$. Then
we have the decomposition ${\rm Sp}=\bigsqcup_{k=0}^{n-1}
{\rm Sp}^{(k)}$.
The element of
${\rm Sp}^{(k)}$ is characterized as follows.

\begin{prop}\label{prop:speccondition}
 An element $\vec{h}=(h_i)\in \{0,1\}^{\bf N}$
belongs to ${\rm Sp}^{(k)}$ if and only if it satisfies
the conditions,
\begin{align}\label{eq:speccond}
\mbox{(i)} &\ \mbox{$h_i+h_{i+1}+
\cdots +h_{i+n-1} \geq 1$ for any $i$.}
\tag{\theequation a}
\\
\mbox{(ii)} &\ %
\mbox{$\vec{h}\approx
\vec{h}^{(k)}$,
where $\vec{h}^{(k)}:=h(\vec{s}^{(k)})=
(\underbrace{0,\dots,0,1}_{k},
(\underbrace{0,\dots,0,1}_{n})^\infty)$.}
\tag{\theequation b}
\end{align}
\addtocounter{equation}{1}
\end{prop}
\par
One can paraphrase the condition (\ref{eq:speccond}a)
as ``There are at most $n-1$ consecutive
$0$'s in $\vec{h}$''.
Here we prove only the necessity of the conditions.
The sufficiency will be proved 
after Prop.\ \ref{prop:onetoone}
in section \ref{sect:correspondence}.

\begin{pf}
Let us assume
$h_i=h_{i+1}=\cdots =
h_{i+n-2}=0$.
Then we have $s_i \prec s_{i+1}\prec\cdots \prec
s_{i+n-1}$, from which $s_{i+n-1}=\epsilon_n$ follows.
Therefore, $h_{i+n-1}=H(s_{i+n-1},s_{i+n})=1$ regardless of
the value of $s_{i+n}$.

The condition (ii) is an immediate consequence of 
(\ref{eq:boundary}).
\end{pf}

Any element $\vec{h}$ of ${\rm Sp}^{(k)}$ is
uniquely written in the form
\begin{equation}
[m_1,\dots,m_r]:=(
\underbrace{0,\dots,0,1}_{m_1},
\dots,
\underbrace{0,\dots,0,1}_{m_r},
(\underbrace{0,\dots,0,1}_{n})^\infty),
\quad
1\leq m_i \leq n,\
m_r\neq n.
\end{equation}
Obviously
\begin{equation*}
{\rm Sp}^{(k)}=
\{ [m_1,\dots,m_r]\mid
r\geq 0,\ 1\leq m_i\leq n,\ m_r\neq n,\ \sum_{i=1}^r m_i
\equiv k \mod n
\}.
\end{equation*}

The surjection $h:{\cal S}^{(k)}
\rightarrow {\rm Sp}^{(k)}$ induces the
decomposition of ${\cal S}^{(k)}$,
\begin{equation}\label{eq:spindecomposition}
{\cal S}^{(k)} = \bigsqcup_{\vec{h}\in
{\rm Sp}^{(k)}} {\cal S}_{\vec{h}},
\quad {\cal S}_{\vec{h}}=
h^{-1}(\vec{h}).
\end{equation}
We call this decomposition the {\it spectral
decomposition} of ${\cal S}^{(k)}$.

Let us introduce the character
of the degeneracy of the spectrum at $\vec{h}$,
\begin{equation}\label{eq:finitecharacter}
\begin{split}
{\rm ch}\,{\cal S}_{\vec{h}}
&=q^{\Delta_k-c/24}
\sum_{\vec{s}\in {\cal S}_{\vec{h}}}
q^{E(\vec{s})}
e^{{\rm wt}(\vec{s})}\\
&=
q^{\Delta_k- c/24 + \sum_{i=1}^\infty
i (h_i-h^{(k)}_i)}
\chi_{\vec{h}},
\quad
\chi_{\vec{h}}=
\sum_{\vec{s}\in {\cal S}_{\vec{h}}}
e^{{\rm wt}(\vec{s})}.
\end{split}
\end{equation}
As is standard in the character theory of $sl_n$,
we regard $\chi_{\vec{h}}=\chi_{\vec{h}}(x)$
as a function of the variables $x_1=e^{\epsilon_1},
x_2=e^{\epsilon_2},\dots,x_n=e^{\epsilon_n}$
with the relation $x_1x_2\cdots x_n=1$.
Due to Theorem \ref{thm:DJKMO},
the character of ${\cal L}(\Lambda_k)$ is
decomposed as
\begin{equation}\label{eq:fullcharacter}
{\rm ch}\, {\cal L}(\Lambda_k)(q,x)=
q^{\Delta_k- c/24}
\sum_{\vec{h}\in {\rm Sp}^{(k)}}
q^{\sum_{i=1}^\infty i (h_i-h^{(k)}_i)} 
\chi_{\vec{h}}(x).
\end{equation}

The main purpose of the paper
is to calculate the characters $\chi_{\vec{h}}$ 
and to show that they are irreducible 
characters of $Y(sl_n)$. To do it,
we make use of a hidden relation 
between spin configurations and skew Young
diagrams.

\section{Skew diagrams and skew Schur functions}

Let us recall the definitions of a skew diagram,
a semi-standard tableau, and
the skew Schur function. We basically
follow the definitions
and notations of \cite{macdonald}.

A {\it partition} $\lambda=(\lambda_1,\lambda_2,\dots,\lambda_m)$
is a non-increasing
sequence of non-negative integers,
$\lambda_1 \geq \lambda_2 \geq \cdots\geq
\lambda_m \geq 0$.
We let $|\lambda|=\sum_{i=1}^m\lambda_i$.
The {\it length $l(\lambda)$ of $\lambda$}
is the number of the non-zero elements in $\lambda$.
As usual, a partition $\lambda$ is represented by a {\it (Young) diagram},
which is denoted by the same symbol $\lambda$.
We conveniently identify the partitions $(\lambda_1,\dots,
\lambda_m)$, $(\lambda_1,\dots,\lambda_m,0)$,
$(\lambda_1,\dots,\lambda_m,0,0)$, {\it etc}.
The {\it conjugate} of a partition $\lambda$
is a partition  $\lambda'$ whose diagram is the
transpose of the diagram of $\lambda$ along the main
diagonal.
For example,
if $\lambda=(4,3,2)$, then its conjugate is $\lambda'=(3,3,2,1)$.

For a pair of partitions $\lambda$ and $\mu$,
we write $\lambda \supset \mu$ if $\lambda_i
\geq \mu_i$ for any $i$.
If $\lambda
\supset \mu$, the diagram $\mu$ is
naturally embedded inside the diagram $\lambda$.
Then the {\it skew (Young) diagram} $\lambda/\mu$
(denoted by $\lambda-\mu$ in \cite{macdonald})
is obtained by subtracting the diagram $\mu$
from the diagram $\lambda$.
For example, if $\lambda=(5,4,4,1)$ and
$\mu=(4,3,2,0)$, then $\lambda/\mu$ looks as follows:
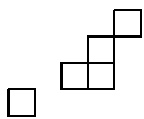
\begin{equation*}
\begin{picture}(50,40)(0,-40)
\put(0,-40){\line(1,0){10}}
\put(0,-30){\line(1,0){10}}
\put(20,-30){\line(1,0){20}}
\put(20,-20){\line(1,0){20}}
\put(30,-10){\line(1,0){20}}
\put(40,0){\line(1,0){10}}
\put(0,-30){\line(0,-1){10}}
\put(10,-30){\line(0,-1){10}}
\put(20,-20){\line(0,-1){10}}
\put(30,-10){\line(0,-1){20}}
\put(40,0){\line(0,-1){30}}
\put(50,0){\line(0,-1){10}}
\end{picture}
\end{equation*}
We set $|\lambda/\mu|=|\lambda|-|\mu|$.

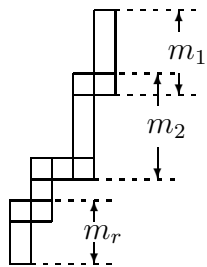
\begin{figure}[bt]
\begin{center}
\setlength{\unitlength}{0.8pt}
\begin{picture}(105,120)(10,-120)
\put(50,0){\line(1,0){10}}
\put(40,-40){\line(1,0){20}}
\put(40,-30){\line(1,0){20}}
\put(20,-70){\line(1,0){30}}
\put(20,-80){\line(1,0){30}}
\put(10,-90){\line(1,0){20}}
\put(10,-100){\line(1,0){20}}
\put(10,-120){\line(1,0){10}}
\put(50,0){\line(0,-1){40}}
\put(60,0){\line(0,-1){40}}
\put(50,-30){\line(0,-1){50}}
\put(40,-30){\line(0,-1){50}}
\put(30,-70){\line(0,-1){30}}
\put(20,-70){\line(0,-1){50}}
\put(10,-90){\line(0,-1){30}}
\multiput(60,0)(6,0){7}{\line(1,0){2}}
\multiput(60,-40)(6,0){7}{\line(1,0){2}}
\multiput(50,-30)(6,0){7}{\line(1,0){2}}
\multiput(50,-80)(6,0){7}{\line(1,0){2}}
\multiput(20,-90)(6,0){7}{\line(1,0){2}}
\multiput(20,-120)(6,0){7}{\line(1,0){2}}
\put(90,-10){\vector(0,1){10}}
\put(90,-30){\vector(0,-1){10}}
\put(85,-23){$m_1$}
\put(80,-45){\vector(0,1){15}}
\put(80,-65){\vector(0,-1){15}}
\put(75,-58){$m_2$}
\put(50,-100){\vector(0,1){10}}
\put(50,-110){\vector(0,-1){10}}
\put(45,-108){$m_r$}
\end{picture}
\end{center}
\caption{A border strip $\langle
m_1,\dots,m_r\rangle$.}
\label{fig:borderstrip}
\end{figure}

We say $\lambda/\mu$ is a skew diagram of
{\it rank} $n$ if
the length of any column of $\lambda/\mu$
does not exceed $n$.
Two boxes in a skew diagram are {\it adjacent}
if they share a common side.
A skew diagram $\lambda/\mu$ is
{\it connected\/} if, for any pair of 
boxes $a$ and $a'$ in $\lambda/\mu$,
 there exits a series
of boxes $b_1=a,b_2,\dots,b_j=a'$ in $\lambda/\mu$
 such that
$b_i$ and $b_{i+1}$ are adjacent.
A skew diagram is called a {\it border strip}
if it is connected, and contains no $2\times 2$
block of boxes.
Let $\langle m_1,\dots,m_r \rangle$ denote
the border strip of $r$ columns such that
the length of $i$th column (from the right)
is $m_i$ (Fig.\ \ref{fig:borderstrip}).

For a skew diagram  $\lambda / \mu$,
we now define the {\it skew Schur function} $s_{\lambda/\mu}$.
In each box of 
a given skew diagram $\lambda/\mu$,
let us inscribe one of the numbers $1,2,\dots,n$.
We call such an arrangement of numbers 
a {\it semi-standard tableau $T$
of shape $\lambda/\mu$},
if it satisfies the following condition:
Let $a$ and $b$ be the inscribed
numbers in any pair of adjacent
boxes. Then,
\begin{align}\label{eq:adjacent}
\mbox{(i)} &\ \mbox{
$a< b$ if $b$ is lower-adjacent to $a$.}
\tag{\theequation a}
\\
\mbox{(ii)} &\ \mbox{
$a\geq b$ if $b$ is left-adjacent to $a$.}
\tag{\theequation b}
\end{align}
\addtocounter{equation}{1}

The $sl_n$-{\it weight} of a semi-standard
tableau $T$ is defined as
\begin{equation}\label{eq:tableaweight}
{\rm wt}'(T)=
\sum_{a=1}^n m_a \cdot \epsilon_a,
\end{equation}
where $m_a$ is the number counting how many
$a$'s are in $T$, and $\epsilon_a$ is given
in (\ref{eq:weights}).
\begin{defn} 
The skew
Schur function $s_{\lambda/\mu}$ is defined as
\begin{equation}\label{eq:skewschurtwo}
s_{\lambda/\mu}(x)=
\sum_{T\in {\rm SST}(\lambda/\mu)} e^{{\rm wt}'(T)},
\quad x_i=e^{\epsilon_i},
\end{equation}
where ${\rm SST}(\lambda/\mu)$ is the set of all the
semi-standard tableaux  of shape $\lambda/ \mu$.
\end{defn}

The following proposition is
well-known.
See section 5 of \cite{macdonald} for a proof.

\begin{prop}\label{prop:equivalent}
The skew Schur function $s_{\lambda / \mu}$ 
is also expressed  as
\begin{equation}\label{eq:skewschur}
s_{\lambda / \mu}(x)=
{\rm det}(e_{\lambda'_i-\mu'_j-i+j}(x))_{%
1\leq i,j \leq r},
\end{equation}
where $r\geq l(\lambda')$,
and  $e_m=e_m(x)$ is the $m${\rm th}
 elementary symmetric
polynomial of variables
$x_1,\dots,x_n$ for
$m=0,\dots,n$, and $e_m=0$ for
other $m$.
\end{prop}

We impose the relation
$x_1 x_2 \cdots x_n=1$ throughout the paper.
Then $e_m$ is 
 the character of the
$m$th fundamental representation
 of $sl_n$ with
the highest weight $\overline\Lambda_m$ for 
$m=1,\dots,n-1$.

The following properties of  $s_{\lambda/
\mu}$ follow  either from (\ref{eq:skewschurtwo})
or from (\ref{eq:skewschur}):
\begin{itemize}
\item[(i)]{ If $\lambda/\mu$ is not a skew diagram of rank $n$, then
$s_{\lambda/\mu}=0$.}
\item[(ii)]{
When  $\mu=(0)$, 
the expression (\ref{eq:skewschur}) reduces to
the Jacobi-Trudi formula
of the ordinary Schur function $s_\lambda$.}
\item[(iii)]{Let $c_{\mu\nu}^\lambda$,
$|\lambda|=|\mu|+|\nu|$, be
the Littlewood-Richardson coefficient, i.e.,
$s_\mu s_\nu=\sum_{\lambda}c_{\mu\nu}^\lambda
s_\lambda$. Then,
%
$
s_{\lambda/\mu}=\sum_\nu c_{\mu\nu}^\lambda
s_\nu.
$
}
\end{itemize}

The {\it conjugate} $s^*_{\lambda/\mu}$ of
the skew Schur function $s_{\lambda/\mu}$
is defined as 
\begin{align}
s^*_{\lambda/\mu}&={\rm det}
(e_{n-\lambda'_i+\mu'_j+i-j})_{%
1\leq i,j \leq r}\\
&=
\sum_{T\in {\rm SST}(\lambda/\mu)}
 e^{-{\rm wt}'(T)}.
\end{align}

It is also possible to express $s^*_{\lambda/\mu}$
as a skew Schur function. 
Suppose $\lambda/\mu$ is of rank $n$
with $\mu=(\mu_1,\dots,\mu_m)$.
Then we have a new pair
\begin{equation*}
\tilde\mu=(\underbrace{\lambda_1,\dots,\lambda_1}_n,
\mu_1,\dots,\mu_m)\supset \lambda,
\end{equation*}
and the {\it compliment of $\lambda/\mu$},
 $(\lambda/\mu)^c:=\tilde\mu/\lambda$,
is also a skew diagram of rank $n$.
The picture below illustrates the example of
$n=4$, $\lambda=(5,4,3,1)$, $\mu=(3,2)$:
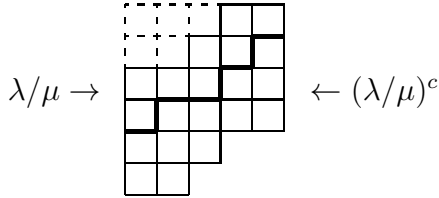
\begin{equation*}
\setlength{\unitlength}{1.2pt}
\begin{picture}(50,60)(0,-60)
\put(30,0){\line(1,0){20}}
\put(20,-10){\line(1,0){10}}
\put(39.8,-9.7){\thicklines\line(1,0){10.4}}
\put(39.8,-10){\thicklines\line(1,0){10.4}}
\put(39.8,-10.3){\thicklines\line(1,0){10.4}}
\put(0,-20){\line(1,0){20}}
\put(29.8,-19.7){\thicklines\line(1,0){10.4}}
\put(29.8,-20){\thicklines\line(1,0){10.4}}
\put(29.8,-20.3){\thicklines\line(1,0){10.4}}
\put(9.8,-29.7){\thicklines\line(1,0){20.4}}
\put(9.8,-30){\thicklines\line(1,0){20.4}}
\put(9.8,-30.3){\thicklines\line(1,0){20.4}}
\put(0,-39.7){\thicklines\line(1,0){10.2}}
\put(0,-40){\thicklines\line(1,0){10.2}}
\put(0,-40.3){\thicklines\line(1,0){10.2}}
\put(30,-40){\line(1,0){20}}
\put(20,-50){\line(1,0){10}}
\put(0,-60){\line(1,0){20}}
\multiput(0,0)(5,0){6}{\line(1,0){1}}
\multiput(4,0)(5,0){6}{\line(1,0){1}}
\multiput(0,-10)(5,0){4}{\line(1,0){1}}
\multiput(4,-10)(5,0){4}{\line(1,0){1}}
\put(30,-10){\line(1,0){10}}
\put(20,-20){\line(1,0){30}}
\put(0,-30){\line(1,0){50}}
\put(0,-40){\line(1,0){30}}
\put(0,-50){\line(1,0){20}}
\put(0,-20){\line(0,-1){40}}
\put(9.7,-30){\thicklines\line(0,-1){10}}
\put(10,-30){\thicklines\line(0,-1){10}}
\put(10.3,-30){\thicklines\line(0,-1){10}}
\put(20,-10){\line(0,-1){10}}
\put(20,-50){\line(0,-1){10}}
\put(30,0){\line(0,-1){10}}
\put(29.7,-20){\thicklines\line(0,-1){10}}
\put(30,-20){\thicklines\line(0,-1){10}}
\put(30.3,-20){\thicklines\line(0,-1){10}}
\put(30,-40){\line(0,-1){10}}
\put(39.7,-10){\thicklines\line(0,-1){10}}
\put(40,-10){\thicklines\line(0,-1){10}}
\put(40.3,-10){\thicklines\line(0,-1){10}}
\put(50,0){\line(0,-1){40}}
\multiput(0,0)(0,-5){4}{\line(0,-1){1}}
\multiput(0,-4)(0,-5){4}{\line(0,-1){1}}
\multiput(10,0)(0,-5){4}{\line(0,-1){1}}
\multiput(10,-4)(0,-5){4}{\line(0,-1){1}}
\put(10,-20){\line(0,-1){40}}
\multiput(20,0)(0,-5){2}{\line(0,-1){1}}
\multiput(20,-4)(0,-5){2}{\line(0,-1){1}}
\put(20,-20){\line(0,-1){30}}
\put(30,-10){\line(0,-1){30}}
\put(40,0){\line(0,-1){40}}
\put(-37,-30){$\lambda/\mu
\rightarrow$}
\put(55,-30){
$\leftarrow(\lambda/\mu)^c$}
\end{picture}
\end{equation*}

\begin{prop}\label{prop:conjugate}
Let $\lambda/\mu$ be a skew 
diagram of  rank $n$.
Then
\begin{equation*}
s^*_{\lambda/\mu}=s_{(\lambda/\mu)^c}.
\end{equation*}
\end{prop}
\begin{pf}
$s_{(\lambda/\mu)^c}=
{\rm det}
(e_{n+\mu'_i-\lambda'_j-i+j})
={\rm det}
(e_{n-\lambda'_i+\mu'_j+i-j})
=s^*_{\lambda/\mu}
$.
\end{pf}

\section{Correspondence between spin configurations
and semi-standard tableaux}\label{sect:correspondence}

We proceed to calculate the character $\chi_{\vec{h}}$
of (\ref{eq:finitecharacter}) using the language of
skew diagrams and  tableaux.

For a given $\vec{h}\in {\rm Sp}$, we associate
a skew diagram $\tilde\kappa(
\vec{h})$ of {\it infinite-size} in the following
procedure:
\begin{itemize}
\item[1.]{Write the first box.}
\item[2.]{Attach the second box under (resp.\
left to) the first box if $h_1=0$ (resp.\ 
$h_1=1$).}
\item[3.]{Similarly
attach the $i+1$th box under (resp.\
left to) the $i$th box if $h_{i}=0$ (resp.\ 
$h_{i}=1$) for $i=2,3,\dots$}
\end{itemize}%
\noindent
Then $\tilde\kappa(\vec{h})$ has
the following properties.
\begin{itemize}
\item[(i)]{It is a border strip.}
\item[(ii)]{It is of rank $n$, i.e., the length of
any column of $\tilde\kappa(\vec{h})$
does not exceed $n$, due to (\ref{eq:speccond}}a).
\item[(iii)]{Due to (\ref{eq:speccond}b),
it has a periodic tail which consists
of length-$n$ columns.}
\end{itemize}
Equivalently, for $\vec{h}=[m_1,\dots,m_r]$,
$\tilde\kappa(\vec{h})$ is the border strip
$\langle m_1,\dots,m_r,n,n,\dots\rangle$.

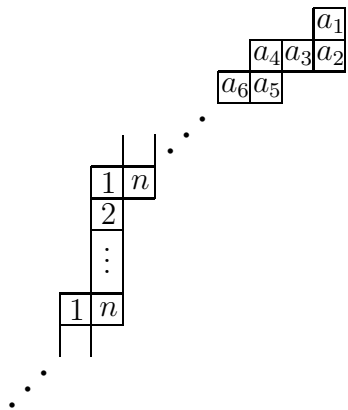
\begin{figure}[bt]
\begin{center}
\setlength{\unitlength}{1.2pt}
\begin{picture}(105,125)(-45,-125)
\put(50,0){\line(1,0){10}}
\put(30,-10){\line(1,0){30}}
\put(20,-20){\line(1,0){40}}
\put(20,-30){\line(1,0){20}}
\put(20,-20){\line(0,-1){10}}
\put(30,-10){\line(0,-1){20}}
\put(40,-10){\line(0,-1){20}}
\put(50,0){\line(0,-1){20}}
\put(60,0){\line(0,-1){20}}
\put(50.5,-10){\makebox(10,10){$a_1$}}
\put(50.5,-20){\makebox(10,10){$a_2$}}
\put(40.5,-20){\makebox(10,10){$a_3$}}
\put(30.5,-20){\makebox(10,10){$a_4$}}
\put(30.5,-30){\makebox(10,10){$a_5$}}
\put(20.5,-30){\makebox(10,10){$a_6$}}
\put(15,-35){\circle*{2}}
\put(10,-40){\circle*{2}}
\put(5,-45){\circle*{2}}
\put(-20,-50){\line(1,0){20}}
\put(-20,-60){\line(1,0){20}}
\put(-20,-70){\line(1,0){10}}
\put(-20,-90){\line(1,0){10}}
\put(-20,-100){\line(1,0){10}}
\put(-30,-90){\line(1,0){10}}
\put(-30,-100){\line(1,0){10}}
\put(0,-40){\line(0,-1){20}}
\put(-10,-40){\line(0,-1){60}}
\put(-20,-50){\line(0,-1){60}}
\put(-30,-90){\line(0,-1){20}}
\put(-9.5,-60){\makebox(10,10){$n$}}
\put(-19.5,-60){\makebox(10,10){$1$}}
\put(-19.5,-70){\makebox(10,10){$2$}}
\put(-20,-82){\makebox(10,10){$\vdots$}}
\put(-19.5,-100){\makebox(10,10){$n$}}
\put(-29.5,-100){\makebox(10,10){$1$}}
\put(-35,-115){\circle*{2}}
\put(-40,-120){\circle*{2}}
\put(-45,-125){\circle*{2}}
\end{picture}
\end{center}
\caption{A semi-standard tableau $\tilde{T}$
 of shape $\tilde\kappa(\vec{h})$
for $\vec{h}=(0, 1, 1, 0, 1,\- \dots, (0, \dots,0,
1)^\infty)$.}
\label{fig:sst}
\end{figure}

A semi-standard
tableau $\tilde T$ of shape $\tilde\kappa(\vec{h})$ is
an arrangement of numbers $1,\dots,n$
in the boxes of $\tilde\kappa(\vec{h})$
obeying the conditions (\ref{eq:adjacent}a) and 
(\ref{eq:adjacent}b),
just in the same way as in the finite-size case.
See Fig.\ \ref{fig:sst}.
Notice that the arrangements in the length-$n$ columns are uniquely
determined, or ``frozen'', because of the semi-standard
condition (\ref{eq:adjacent}a).
Since $\tilde\kappa(\vec{h})$ is a border strip,
one can give a  total
ordering of the boxes in it
from the right to the left and from the top to the
bottom in the unique way.
Let
\begin{equation*}
(a_i)=(a_1,a_2,\dots,(1,2,\dots,n)^\infty).
\end{equation*}
be the sequence of  the content of $\tilde{T}$
along the total ordering of the boxes.
Now we have a natural map $\tilde\varphi$ from the set
of the semi-standard tableaux of shape
$\tilde\kappa(\vec{h})$ to the set of the spin
configurations ${\cal S}$ defined by
\begin{equation}\label{eq:bijectionphi}
\tilde\varphi(\tilde T):=(\epsilon_{a_i})
=(\epsilon_{a_1},\epsilon_{a_2},\dots,
(\epsilon_1,\epsilon_2,\dots,\epsilon_n)^\infty).
\end{equation}

A key observation of this paper is
\begin{prop}\label{prop:onetoone}
The map $\tilde\varphi$ gives a one-to-one correspondence between
the semi-stand\-ard
 tableaux of shape $\tilde\kappa(\vec{h})$ and the spin configurations
in ${\cal S}_{\vec{h}}$. 
\end{prop}

\begin{pf}
A necessary and sufficient condition for
a sequence
$(a_1,a_2,\dots)$, $a_i\in \{
1,\dots,n\}$ to be the sequence of
the content of
a semi-standard tableau of shape
$\tilde\kappa(\vec{h})$ is
\begin{equation}\label{eq:tableaucondition}
a_{i+1} > a_i,\
({\rm resp.\ }a_{i+1}\leq a_i)
\Longleftrightarrow
h_i=0,\
({\rm resp.\ } h_i=1).
\end{equation}
This is an immediate consequence of the construction
of $\tilde\kappa(\vec{h})$ and (\ref{eq:adjacent}a,b).
On the other hand (\ref{eq:tableaucondition})
is also a necessary and sufficient condition
for a sequence
$(\epsilon_{a_1},\epsilon_{a_2},\dots)$
to belong to ${\cal S}_{\vec{h}}$ because of
(\ref{eq:localenergy}), (\ref{eq:maph}),
and (\ref{eq:spindecomposition}). 
\end{pf}

The first application of Prop.\ \ref{prop:onetoone} is
to prove the sufficiency
part of Prop.\ \ref{prop:speccondition}.
\begin{pf}
By Prop.\ \ref{prop:onetoone} we have only to
show that there exists at least one
semi-standard tableau of shape
$\tilde\kappa(\vec{h})$ for any $\vec{h}$.
In fact, for a given $\vec{h}$,
a semi-standard tableau of shape
$\tilde\kappa(\vec{h})$ is obtained by
filling the boxes by $1,2,3,\dots$
from the top to the bottom in each column.
This completes the proof of Prop.\ \ref{prop:speccondition}.
\end{pf}

For our purpose it is convenient to define a
``finite part'' $\kappa(\vec{h})$ 
of the infinite diagram
$\tilde\kappa(\vec{h})$ by cutting off
its periodic tail.
Namely, for $\vec{h}=[m_1,\dots,m_r]$
we define $\kappa(\vec{h})=\langle
m_1,\dots,m_r\rangle$.
(For $\vec{h}=\vec{h}^{(0)}$, 
$\kappa(\vec{h}^{(0)})$ is the empty diagram $\emptyset$.)
It is clear that the map $\kappa:
\vec{h}\mapsto \kappa(\vec{h})$
is injective, and we get the following
description
of the space of the spectrum ${\rm Sp}$
in terms of border strips.
\begin{thm}
The space ${\rm Sp}$ is parametrized by
the border strips of rank $n$ such that
the lengths of their leftmost columns
are less than $n$.
\end{thm}

The following lemma is obvious.

\begin{lem}\label{lemma:bijection}
There is a one-to-one correspondence
between the semi-standard tabl\-eaux of
shape $\tilde\kappa(\vec{h})$ and the ones
of shape $\kappa(\vec{h})$. The correspondence
is given by the restriction of a
semi-standard tableau of shape $\tilde\kappa
(\vec{h})$ on $\kappa(\vec{h})$.
\end{lem}

Combining the bijection in
 Lemma \ref{lemma:bijection} with the bijection
$\tilde\varphi$ of (\ref{eq:bijectionphi}),
we obtain a bijection
\begin{equation*}
\varphi: 
{\rm SST}(\kappa(\vec{h}))
\rightarrow 
{\cal S}_{\vec{h}}.
\end{equation*}

\begin{prop}\label{prop:weightpreserve}
The bijection $\varphi:
{\rm SST}(\kappa(\vec{h}))
\rightarrow 
{\cal S}_{\vec{h}}$ is weight-preserving,
i.e., for any $T\in {\rm SST}(\kappa(\vec{h}))$,
${\rm wt}(\varphi(T))={\rm wt}'(T)$ holds.
\end{prop}
\begin{pf}
Let $(a_1,\dots,a_m)$ be the content of $T$
aligned along our total order of the boxes
in $\kappa(\vec{h})$.
Then $\varphi(T)=
(\epsilon_{a_1},\dots,\epsilon_{a_m},
(\epsilon_1,\dots,\epsilon_n)^\infty)$.
{}From Prop.\ \ref{prop:weight} and (\ref{eq:tableaweight}),
we have ${\rm wt}(\varphi(T))=\sum_{i=1}^m\epsilon_{a_i}
={\rm wt}'(T)$.
\end{pf}

Now we state the first half of our main theorem.
\begin{thm}
(i) The character $\chi_{\vec{h}}$ of
${\cal S}_{\vec{h}}$ is equal to the skew
Schur function $s_{\kappa(\vec{h})}$.
\hfill\break
(ii) Let $\vec{h}=[m_1,\dots,m_r]
\in {\rm Sp}$. Then
\begin{equation}\label{eq:schurformula}
s_{\kappa(\vec{h})}
=s_{\langle m_1,\dots,m_r\rangle}
=
\left|
\begin{array}{llllll}
e_{m_r} & e_{m_r+m_{r-1}}&
&\cdots&& e_{m_r+\cdots+m_1}\\
1 & e_{m_{r-1}}\\
0 & 1\\
&0&&&&\vdots\\
&&\ddots&\ddots\\
&&0&1& e_{m_2}&e_{m_2+m_1}\\ 
&&&0&1& e_{m_1} 
\end{array}
\right|\,.
\end{equation}
\hfill\break
(iii)
The character of the level 1 integrable module
${\cal L}(\Lambda_k)$ of $\widehat{sl}_n$
decomposes as
\begin{align}\label{eq:characterdecomposition}
{\rm ch}\, {\cal L}(\Lambda_k)(q,x)
&=
q^{-{1\over24}c}
\sum_{\kappa\in BS\atop |\kappa|\equiv k\, {\rm mod}\, n}
q^{{1\over2n}|\kappa|(n-|\kappa|)+t(\kappa)} 
s_{\kappa}(x)
\tag{\theequation a}\\
&=
q^{-{1\over24}c}
\sum_{\kappa\in BS\atop |\kappa|\equiv n-k \, {\rm mod}\, n}
q^{{1\over2n}|\kappa|(n-|\kappa|)+t(\kappa)} 
s_{\kappa^c}(x)
\tag{\theequation b}\,,
\end{align}
where $BS$ is the set of all the border strips
$\kappa=\langle m_1,\dots,m_r\rangle$ of rank $n$
with $m_r<n$, and $t(\kappa)=\sum_{i=1}^{r-1}
(r-i)m_i$.
\end{thm}
\begin{pf}
The property (i) is an immediate consequence
of Prop.\ \ref{prop:weightpreserve}.
To show (ii), notice that 
the skew diagram $\langle m_1,\dots,m_r\rangle$ is
represented as $\lambda/\mu$ with
a pair $\lambda \supset \mu$
such that
\begin{equation*}
\lambda'_i
=m_1+\cdots+
m_{r+1-i}-r+i, 
\quad
\mu'_i
=m_1+\cdots+
m_{r-i}-r+i.
\end{equation*}
Substituting them into 
(\ref{eq:skewschur}) we obtain the formula (\ref{eq:schurformula}).
The property (iii) follows from (\ref{eq:fullcharacter}),
Prop.\ \ref{prop:conjugate}, the property (i),
and the fact that for
$\vec{h}=[m_1,\dots,m_r]\in {\rm Sp}^{(k)}$
\begin{equation*}\label{eq:plainformula}
\Delta_k + \sum_{i=1}^\infty
i(h_i-h^{(k)}_i)
={1\over2n}m(n-m) + \sum_{i=1}^{r-1}(r-i) m_i,
\quad m=\sum_{i=1}^r m_i.
\end{equation*}
\end{pf}
\addtocounter{equation}{1}

Two expressions
(\ref{eq:characterdecomposition}a)
and (\ref{eq:characterdecomposition}b)
differ from each other
when $n\geq 3$.

\section{Yangian characters}

In this section we show, based on the
result of \cite{nazarov}, that the characters
$\chi_{\vec{h}}=s_{\kappa(\vec{h})}$ are
irreducible characters of the Yangian
algebras $Y(gl_n)$ and $Y(sl_n)$.

The Yangian of $gl_n$, $Y(gl_n)$,
is an algebra generated by $t_{ij}(r)$,
$i,j=1,\dots,n$, $r\in {\bf Z}_{\geq 0}$
with the relations
\begin{equation*}
[t_{ij}(r),t_{kl}(s-1)]
-[t_{ij}(r-1),t_{kl}(s)] =t_{kj}(s-1)t_{il}(r-1)
-t_{kj}(r-1)t_{il}(s-1),
\end{equation*}
where $t_{ij}(-1)=\delta_{ij}1$,
$t_{ij}(-2)=0$. The elements
$t_{ij}(0)$ generate
the universal enveloping
algebra of $gl_n$.
 Therefore $sl_n$ acts on  $Y(gl_n)$-modules.

Consider a  pair of partitions  $\lambda \supset
\mu$ with  $\lambda=(\lambda_1,\dots,
\lambda_{N+n})$, $\mu=(\mu_1,\dots,
\mu_N)$, $N\geq 1$.
Let $V_\lambda$ be the irreducible
$gl_{N+n}$-module associated to $\lambda$,
$V_\mu$ be the irreducible
$gl_{N}$-module associated to $\mu$,
and $V_{\lambda,\mu}$
be the space of the multiplicity of 
$V_\mu$ in $V_\lambda$
under the standard embedding $gl_N
\subset gl_{N+n}$.
There is an
irreducible action of $Y(gl_n)$ on
the space $V_{\lambda,\mu}$,
having  a remarkable property
(A module with such  property
is called a {\it tame module} \cite{nazarov}):
\begin{prop}[\cite{cherednik,nazarov}]
\label{prop:gz}
A maximal commutative subalgebra of 
$Y(gl_n)$, called the Gelfand-Zetlin
(GZ) algebra, acts on $V_{\lambda,\mu}$
in a semi-simple way. 
Furthermore, a basis diagonalizing the GZ algebra
is labeled by the {\it GZ schemes} of $V_{\lambda,\mu}$.
\end{prop}

A GZ scheme $\backslash \Lambda /$
 of $V_{\lambda,\mu}$ is an
array of integers $\lambda_{mi}$,
\begin{equation*}
\backslash \Lambda / = 
\setlength{\unitlength}{1.2pt}
\lower22pt
\hbox{
\begin{picture}(140,50)(0,-45)
\put(0,0){$\lambda_{n1}$}
\put(20,0){$\lambda_{n2}$}
\put(53,0){$\cdots$}
\put(100,0){$\lambda_{n,N+n}$}
\put(10,-15){$\lambda_{n-1,1}$}
\put(53,-15){$\cdots$}
\put(90,-15){$\lambda_{n-1,N+n-1}$}
\put(53,-30){$\cdots$}
\put(30,-45){$\lambda_{01}$}
\put(53,-45){$\cdots$}
\put(70,-45){$\lambda_{0N}$}
\end{picture}
}
=
\lower22pt
\hbox{
\begin{picture}(140,50)(0,-45)
\put(0,0){$\lambda_{1}$}
\put(20,0){$\lambda_{2}$}
\put(53,0){$\cdots$}
\put(100,0){$\lambda_{N+n}$}
\put(10,-15){$\lambda_{n-1,1}$}
\put(53,-15){$\cdots$}
\put(90,-15){$\lambda_{n-1,N+n-1}$}
\put(53,-30){$\cdots$}
\put(30,-45){$\mu_{1}$}
\put(53,-45){$\cdots$}
\put(70,-45){$\mu_{N}$}
\end{picture}
}
\end{equation*}
satisfying the condition
$  \lambda_{mi}\geq
\lambda_{m-1,i}\geq
\lambda_{m,i+1}$.
The $sl_n$-weight of the basis vector labeled
by $\backslash \Lambda /$ is 
\begin{equation*}
\sum_{m=1}^{n}
\left\{
\sum_{i=1}^{N+m}
\lambda_{mi}
-
\sum_{i=1}^{N+m-1}
\lambda_{m-1,i}
\right\}\epsilon_m.
\end{equation*}
\begin{lem}\label{lemma:GZ}
There is a  weight-preserving,
one-to-one correspondence
between the GZ schemes of $V_{\lambda,\mu}$
and the semi-standard tableaux of shape 
$\lambda/\mu$.
\end{lem}

The correspondence is described as follows.
For a given GZ scheme $\backslash\Lambda /$,
we have a sequence of partitions,
\begin{equation*}
\lambda^{(0)}=\mu
\subset \lambda^{(1)}
\subset \lambda^{(2)}
\subset \cdots
\subset \lambda^{(n)}=\lambda,
\quad
\lambda^{(m)}
=(\lambda_{m1},\lambda_{m2},
\dots,\lambda_{m,N+m}).
\end{equation*}
A semi-standard tableau of shape
$\lambda/\mu$ is obtained by
inscribing the number $m$
in $\lambda^{(m)}/\lambda^{(m-1)}$,
which is a part of the diagram
$\lambda/\mu$. For example,
\begin{equation*}
\setlength{\unitlength}{1.2pt}
\lower25pt
\hbox{
\begin{picture}(105,50)(0,-45)
\put(0,0){$5$}
\put(20,0){$4$}
\put(40,0){$4$}
\put(60,0){$1$}
\put(80,0){$0$}
\put(100,0){$0$}
\put(10,-15){$5$}
\put(30,-15){$4$}
\put(50,-15){$4$}
\put(70,-15){$0$}
\put(90,-15){$0$}
\put(20,-30){$4$}
\put(40,-30){$4$}
\put(60,-30){$2$}
\put(80,-30){$0$}
\put(30,-45){$4$}
\put(50,-45){$3$}
\put(70,-45){$2$}
\end{picture}
}
\mapsto
\quad
\setlength{\unitlength}{1.2pt}
\lower20pt
\hbox{
\begin{picture}(50,40)(0,-40)
\put(40,0){\line(1,0){10}}
\multiput(0,0)(5,0){8}{\line(1,0){1}}
\multiput(4,0)(5,0){8}{\line(1,0){1}}
\put(30,-10){\line(1,0){20}}
\multiput(0,-10)(5,0){6}{\line(1,0){1}}
\multiput(4,-10)(5,0){6}{\line(1,0){1}}
\put(20,-20){\line(1,0){20}}
\multiput(0,-20)(5,0){4}{\line(1,0){1}}
\multiput(4,-20)(5,0){4}{\line(1,0){1}}
\put(0,-30){\line(1,0){10}}
\put(20,-30){\line(1,0){20}}
\multiput(10,-30)(5,0){2}{\line(1,0){1}}
\multiput(14,-30)(5,0){2}{\line(1,0){1}}
\put(0,-40){\line(1,0){10}}
\put(0,-30){\line(0,-1){10}}
\multiput(0,0)(0,-5){6}{\line(0,-1){1}}
\multiput(0,-4)(0,-5){6}{\line(0,-1){1}}
\put(10,-30){\line(0,-1){10}}
\multiput(10,0)(0,-5){6}{\line(0,-1){1}}
\multiput(10,-4)(0,-5){6}{\line(0,-1){1}}
\put(20,-20){\line(0,-1){10}}
\multiput(20,0)(0,-5){4}{\line(0,-1){1}}
\multiput(20,-4)(0,-5){4}{\line(0,-1){1}}
\put(30,-10){\line(0,-1){20}}
\multiput(30,0)(0,-5){2}{\line(0,-1){1}}
\multiput(30,-4)(0,-5){2}{\line(0,-1){1}}
\put(40,0){\line(0,-1){30}}
\put(50,0){\line(0,-1){10}}
\put(40,-10){\makebox(10,10){$2$}}
\put(30,-20){\makebox(10,10){$1$}}
\put(30,-30){\makebox(10,10){$2$}}
\put(20,-30){\makebox(10,10){$2$}}
\put(0,-40){\makebox(10,10){$3$}}
\end{picture}
}
\end{equation*}
It is easy to check that this map is
bijective and weight-preserving.

It follows from Prop.\ \ref{prop:gz}
and Lemma \ref{lemma:GZ}
that the $sl_n$-character of $V_{\lambda,\mu}$
is equal to $s_{\lambda/\mu}$. In particular,
we have
\begin{thm} The character $\chi_{\vec{h}}
=s_{\kappa(\vec{h})}$ is the $sl_n$-character
of the irreducible $Y(gl_n)$-module $V_{\lambda,
\mu}$ with $\lambda/\mu=\kappa(\vec{h})$.
\end{thm}

The Yangian of $sl_n$, $Y(sl_n)$, is generated by
$x^\pm_{ik}$, $h_{ik}$, $i=1,\dots,n-1$, $k\in 
{\bf Z}_{\geq 0}$ with the relations
\begin{equation*}
\begin{split}
&[h_{ik},h_{jl}]=0,\quad
[h_{i0},x^\pm_{jl}]=\pm A_{ij}
x^\pm_{jl},\quad
[x^+_{ik},x^-_{jl}]=\delta_{ij}h_{ik+l},\\
&[h_{ik+1},x^\pm_{jl}]-[h_{ik},x^\pm_{jl+1}]=
\pm{1\over2}A_{ij}
(h_{ik}x^\pm_{jl} + x^\pm_{jl}h_{ik}),\\
&[x^\pm_{ik+1},x^\pm_{jl}]-[x^\pm_{ik},x^\pm_{jl+1}]=
\pm{1\over2}A_{ij}
(x^\pm_{ik}x^\pm_{jl} + x^\pm_{jl}x^\pm_{ik}),\\
&\sum_{\sigma: {\rm permutation}}
\left[x^\pm_{ik_{\sigma(1)}},
\left[x^\pm_{ik_{\sigma(2)}},
\dots,
\left[x^\pm_{ik_{\sigma(1-A_{ij})}},
x^\pm_{jl}\right]\dots\right]
\right]=0,
\end{split}
\end{equation*}
where $A_{ij}$ is the Cartan matrix of $sl_n$.
An irreducible finite-dimensional module
of $Y(sl_n)$ is characterized by $n-1$ monic polynomials
(the {\it Drinfel'd polynomials}), $P_1(u),\dots,P_{n-1}(u)$.
The polynomial $P_i(u)$ describes the action of 
$h_i(u)=1+\sum_{k=0}^\infty h_{ik}u^{-k-1}$ on a
highest weight vector $v$  as
$h_i(u)v=(P_i(u+1)/ P_i(u))v$ \cite{drinfeld}.

\begin{prop}[\cite{nazarov}]\label{prop:nazarov}
There is a 
one-parameter family of irreducible
$Y(sl_n)$-module structures on $V_{\lambda,\mu}$
with a parameter $b\in {\bf C}$, whose
Drinfel'd polynomials are  \cite{nazarov}\footnote{%
The Drinfel'd polynomials here
are the  one in \cite{drinfeld}. The
convention in \cite{nazarov} is slightly
different. Their $P_i(u)$ is
equal to $(-1)^{{\rm deg}P_i}
P_i(-u-n/4+i/2)$ here.
}
\begin{equation}\label{eq:drinfeld}
P_i(u)=\prod_{j=1\atop\lambda'_j
-\mu'_j=i}^{\lambda_1}
\left(u+{1\over2}(\lambda'_j
+\mu'_j)-j+{1\over2}+b\right).
\end{equation}
\end{prop}

\begin{figure}
\begin{center}
\setlength{\unitlength}{1.2pt}
\begin{picture}(100,130)(-40,-65)
\put(40,0){\line(1,0){10}}
\multiput(0,0)(5,0){8}{\line(1,0){1}}
\multiput(4,0)(5,0){8}{\line(1,0){1}}
\put(30,-10){\line(1,0){20}}
\multiput(0,-10)(5,0){6}{\line(1,0){1}}
\multiput(4,-10)(5,0){6}{\line(1,0){1}}
\put(20,-20){\line(1,0){20}}
\multiput(0,-20)(5,0){4}{\line(1,0){1}}
\multiput(4,-20)(5,0){4}{\line(1,0){1}}
\put(0,-30){\line(1,0){10}}
\put(20,-30){\line(1,0){20}}
\multiput(10,-30)(5,0){2}{\line(1,0){1}}
\multiput(14,-30)(5,0){2}{\line(1,0){1}}
\put(0,-40){\line(1,0){10}}
\put(-30,0){\vector(1,0){90}}
\put(-20,65){\line(0,-1){115}}
\put(0,-50){\vector(0,1){115}}
\put(0,-30){\line(0,-1){10}}
\multiput(0,0)(0,-5){6}{\line(0,-1){1}}
\multiput(0,-4)(0,-5){6}{\line(0,-1){1}}
\put(10,-30){\line(0,-1){10}}
\multiput(10,0)(0,-5){6}{\line(0,-1){1}}
\multiput(10,-4)(0,-5){6}{\line(0,-1){1}}
\put(20,-20){\line(0,-1){10}}
\multiput(20,0)(0,-5){4}{\line(0,-1){1}}
\multiput(20,-4)(0,-5){4}{\line(0,-1){1}}
\put(30,-10){\line(0,-1){20}}
\multiput(30,0)(0,-5){2}{\line(0,-1){1}}
\multiput(30,-4)(0,-5){2}{\line(0,-1){1}}
\put(40,0){\line(0,-1){30}}
\put(50,0){\line(0,-1){10}}
\put(-20,-10){\line(1,-1){25}}
\put(-20,20){\line(1,-1){45}}
\put(-20,35){\line(1,-1){55}}
\put(-20,60){\line(1,-1){65}}
\put(40,-10){\makebox(10,10){$\bullet$}}
\put(30,-25){\makebox(10,10){$\bullet$}}
\put(20,-30){\makebox(10,10){$\bullet$}}
\put(0,-40){\makebox(10,10){$\bullet$}}
\put(-25,55){\makebox(10,10){$\bullet$}}
\put(-25,30){\makebox(10,10){$\bullet$}}
\put(-25,15){\makebox(10,10){$\bullet$}}
\put(-25,-15){\makebox(10,10){$\bullet$}}
\put(-35,-60){$\textstyle x=b$}
\put(-60,-12){$-3-b$}
\put(-60,18){$\phantom{-2}-b$}
\put(-60,32){$\phantom{-}{3\over2}-b$}
\put(-60,58){$\phantom{-}4-b$}
\end{picture}
\end{center}
\caption{An example: $\lambda=(5,4,4,1)$,
$\mu=(4,3,2,0)$. The Drinfel'd polynomials
are $P_1(u)=(u+3+b)
(u+b) (u-4+b)$,
$P_2(u)=u-{3\over2}+b$, and $P_i(u)=1$ for $i\geq 3$.
The zeros of $P_i(u)$ are identified with 
the intersections of the line $x=b$
and the diagonal lines passing through the
middle points of the columns of length $i$.}
\label{fig:drinfeld}
\end{figure}
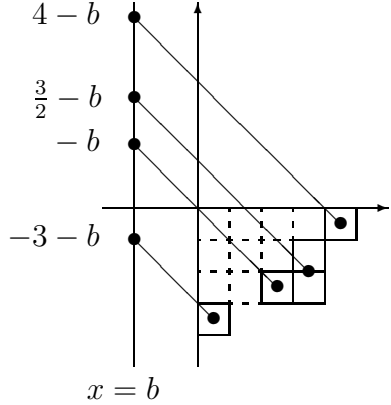

These $Y(sl_n)$-module structures on $V_{\lambda,\mu}$
 are the ones induced by
a one-parameter family of embeddings of $Y(sl_n)$
into $Y(gl_n)$.
There is a simple pictorial interpretation
of the zeros of $P_i(u)$  as shown
in Fig.\ \ref{fig:drinfeld}.

As a corollary of Prop.\ \ref{prop:nazarov},
we obtain the second half of our main theorem:
\begin{thm}\label{thm:main}
 The character $\chi_{\vec{h}}=
s_{\kappa(\vec{h})}$ is the $sl_n$-character
of the one-parameter family of
the irreducible $Y(sl_n)$-modules
whose Drinfel'd polynomials are given by
(\ref{eq:drinfeld}) with
$\lambda/\mu=\kappa(\vec{h})$.
\end{thm}

We have shown that  the characters $\chi_{\vec{h}}$
of the degeneracy of the spectrum are irreducible
$Y(sl_n)$ characters.
Furthermore, the $Y(sl_n)$-module structure of ${\cal L}(\Lambda_k)$
partially
studied in  \cite{schoutens} agrees with
(\ref{eq:characterdecomposition}a).
Based on these strong evidences, we conjecture that
\begin{conj}
The decomposition (\ref{eq:characterdecomposition}a)
describes the $Y(sl_n)$-module structure
on ${\cal L}(\Lambda_k)$ of
\cite{schoutens}.
\end{conj}

\section{The relation with
the spectrum in other spin models}

\subsection{The Haldane-Shastry model}
The $sl_n$ Haldane-Shastry (HS) model is a lattice model
with the Hamiltonian
\begin{equation*}
{\cal H}_{\rm HS}=
\sum_{1\leq j \neq k \leq N}
{x_j x_k \over (x_j-x_k)(x_k-x_j)}(P_{jk}-1),
\quad x_j=e^{2\pi \sqrt{-1} j/N}
\end{equation*}
acting on the vector space $V^{\otimes N}$,
$V={\bf C}^n$, where $P_{jk}$  exchanges
 the $j$th and $k$th components of $V^{\otimes N}$.
There is an action of $Y(sl_n)$ on $V^{\otimes N}$
which commutes with the Hamiltonian ${\cal H}_{\rm HS}$
\cite{haldanePRL}. By this action the space $V^{\otimes
N}$ decomposes as \cite{bernardJP}
\begin{equation}\label{eq:decompose}
V^{\otimes N}\simeq
\bigoplus_{d\in M_N} W_d,
\end{equation}
where $W_d$ are certain irreducible $Y(sl_n)$ modules
described below.
$M_N$ is the set of the binary sequences
$d=(d_1,\dots,d_{N-1})$, $d_i\in \{0,1\}$,
such that there are at most $n-1$ consecutive
$1$'s.
The eigenvalue of ${\cal H}_{\rm HS}$ on the eigenspace
$W_d$ is  given by
$\sum_{i=1}^{N-1} i d_i (id_i-N)$.

We notice that the condition for $d\in M_N$ turns 
into (\ref{eq:speccond}a) through the
identification $h_i=1-d_i$.\ \footnote{This
 intriguing relation between the spectrum of
the HS model and the vertex model was first indicated
by Bernard \cite{bernardtalk} in the $sl_2$ case.}
For a given $d\in M_N$, let
\begin{equation*}
\vec{h}_d=(1-d_1,1-d_2,
\dots,1-d_{N-1},1,(\underbrace{0,\dots,0,1}_{n})^\infty)
\in {\rm Sp}.
\end{equation*}
We translate the description of the
$Y(sl_n)$-module structure of $W_d$ in \cite{bernardJP}
into our language as follows: 
\begin{prop}[\cite{bernardJP}]\label{prop:bernard}
As a $Y(sl_n)$-module, the eigenspace $W_d$ is isomorphic
to the irreducible module whose Drinfel'd polynomials are
given by (\ref{eq:drinfeld}),
where in (\ref{eq:drinfeld}) 
$\lambda_1$ and
$b$ are certain 
constants independent of $d$,
and $\lambda/\mu=\kappa(\vec{h}_d)$.
\end{prop}

Comparing Prop.\ \ref{prop:bernard} with Prop.\
\ref{prop:nazarov} and Theorem \ref{thm:main},
we obtain the following proposition,
which answers the questions of the $sl_n$-module content
of $W_d$ and its factorizability
asked in \cite{bernardJP,
ha,haldanePRL,hikami}.
\begin{prop}\label{cor:motif}
Let $\vec{h}_d=[m_1,\dots,m_r]$. Then\hfill\break
(i)  ${\rm ch}\, W_d = s_{\langle m_1,\dots,m_r\rangle}$.
\hfill\break
(ii) $s_{\langle m_1,\dots,m_r\rangle}
=s_{\langle m_1,\dots,m_i\rangle}
 s_{\langle m_{i+1},\dots,m_r\rangle}$
if $m_i+m_{i+1}\geq n+1$.
\end{prop}
\begin{pf}
We only need to prove the property (ii), which 
follows from (\ref{eq:schurformula}).
\end{pf}

\subsection{The Polychronakos model}
There is another relevant spin model,
called the Polychronakos  model.
The $sl_n$ Polychronakos model has the
Hamiltonian
\begin{align}\label{eq:PFhamiltonian}
{\cal H}_{\rm P}&=
\sum_{1\leq j < k \leq N}
{1\over (x_j-x_k)^2}(P_{jk}-1)+E_N,\\
E_N&={n-1\over 2n}N^2-{\overline{N}(n-\overline{N})
\over 2n},\end{align}
acting on
$V^{\otimes N}$, $V={\bf C}^n$,
where $x_1,\dots,x_N$ are the zeros of
the Hermite polynomial of degree $N$,
and 
$\overline{N}\equiv N \mod n,\
0\leq \overline{N} \leq n-1$.
The constant $E_N$ is added to make
the ground state energy zero.
Define the partition function of the 
Polychronakos model
\begin{equation}
Z_N^{\rm P}(q,x)=
{\rm tr}_{V^{\otimes N}}\, q^{{\cal H}_{\rm P}}
\prod_{i=1}^{n-1}x_i^{h_i+\cdots
+h_{n-1}},
\end{equation}
where $h_i$'s are the standard basis of the
Cartan subalgebra of $sl_n$.
It  is shown in \cite{polychronakos} that
\begin{equation}
Z_N^{\rm P}(q,x)=
q^{E_N}H_N(q^{-1},x),
\end{equation}
where $H_N(q,x)$ is a generalization of the {\it
Rogers-Szeg\"o polynomial} \cite{andrews},
\begin{gather}
H_N(q,x)=\sum_{
k_i\in {\bf Z}_{\geq 0}
\atop
k_1+\cdots + k_n=N
}
{(q)_N \over (q)_{k_1}(q)_{k_2}\cdots
(q)_{k_n}}
x_1^{k_1}\cdots x_n^{k_n},\quad
(q)_k=\prod_{i=1}^k (1-q^i).
\end{gather}

Again there is an action of $Y(sl_n)$
on $V^{\otimes N}$, which
commutes with  ${\cal H}_{\rm P}$
\cite{hikami}. Based on a numerical
study, it was conjectured in \cite{hikami}
that
\begin{itemize}
\item[1.]{\it
As a $Y(sl_n)$-module, the spin space $V^{\otimes N}$
decomposes exactly in the same way as
in (\ref{eq:decompose}).}
\item[2.]{\it
The eigenvalue $E_d$ of
${\cal H}_{\rm P}$ on $W_d$ is 
\begin{equation}\label{eq:PFenergy}
E_d=-\sum_{i=1}^{N-1}i d_i + E_N.
\end{equation}
}
\end{itemize}
Let us introduce the sets,
\begin{align}
{\rm Sp}_N&=\{ \vec{h} \in {\rm Sp}^{(\overline{N})}
\mid \mbox{$h_i=h_i^{(\overline{N})}$ for $i \geq N\}$},\\
{\cal S}_N&=\{ \vec{s} \in {\cal S}^{(\overline{N})}
\mid \mbox{$s_i=s_i^{(\overline{N})}$ for $i \geq N+1\}$}.
\end{align}

\begin{lem}\label{lem:PFvertex}
Let $N$ be a positive integer.
Then
\begin{itemize}
\item[\it (i)]{The map $d\mapsto \vec{h}_d$ is
a bijection from $M_N$ to ${\rm Sp}_N$.}
\item[\it (ii)]{$h^{-1}({\rm Sp}_N)
={\cal S}_N$.}
\item[\it (iii)]{For any $d\in M_N$ and $\vec{s}\in {\cal S}_{\vec{h}_d}$,
$E_d=E(\vec{s})$.} 
\item[\it (iv)]{The $sl_n$-character of $W_d$ is
$s_{\kappa(\vec{h}_d)}=\chi_{\vec{h}_d}$.}
\end{itemize}
\end{lem}

We define
\begin{equation}\label{eq:finitepartition}
Z_{N}^{\rm vertex}(q,x)=
\sum_{\vec{s}\in {\cal S}_N}
q^{E(\vec{s})} e^{{\rm wt}(\vec{s})},
\end{equation}
Then Lemma \ref{lem:PFvertex} and the above conjecture of
\cite{hikami} claims the identification,
\begin{equation}\label{eq:PFvertex}
Z_N^{\rm P}=Z_N^{\rm vertex}.
\end{equation}
This exact equivalence of two spectra is
intriguing, because
the Hamiltonian for $Z_N^{\rm vertex}$
is of nearest-neighborhood type,
while ${\cal H}_P$ is not.
In fact, the following theorem proves
 (\ref{eq:PFvertex}) directly,
thereby providing a further
support for the conjecture of \cite{hikami}.

\begin{thm}\label{thm:rogers}
For any nonnegative integer $N$,
we have
\begin{equation*}\label{thm:sumrule}
\sum_{\vec{h}\in {\rm Sp}_N}
q^{\sum_{i=1}^\infty i (h_i-h_i^{(\overline{N})})}
s_{\kappa(\vec{h})}(x)
=q^{E_N} H_N(q^{-1},x).
\end{equation*}
\end{thm}
A proof of Theorem \ref{thm:rogers}
is given in Appendix A.

As a corollary of Theorem \ref{thm:rogers}
we get an expression
\begin{equation}\label{eq:limitformula}
q^{\Delta_k - c/24}
\sum_{\vec{s}\in {\cal S}^{(k)}}
 q^{E(\vec{s})}
e^{{\rm wt}(\vec{s})}
=
\lim_{N\rightarrow \infty\atop
N\equiv k\ {\rm mod}\ n}
q^{\Delta_k-c/24+E_N} H_N(q^{-1},x).
\end{equation}
The right hand side of (\ref{eq:limitformula}) converges to
(cf. \cite{hikami,polychronakos})
\begin{equation*}
{1\over \eta(q)^{n-1}}
\sum_{k_i\in {\bf Z}-{k\over n}\atop
k_1+\cdots +k_n=0}
q^{{1\over2} k_1^2+\cdots +{1\over2}k_n^2}
x_1^{k_1}\cdots x_n^{k_n},
\quad \eta(q)=q^{1\over24}(q)_\infty,
\end{equation*}
which is equal to
\begin{equation}\label{eq:kacfrenkel}
{1\over \eta(q)^{n-1}}
\sum_{\gamma\in \oplus_{i=1}^{n-1}
{\bf Z}\alpha_i}
q^{{1\over2}|\overline{\Lambda}_k+\gamma|^2}
e^{\overline{\Lambda}_k+\gamma},
\end{equation}
where $\alpha_i$'s, $|\alpha_i|^2=2$,
are the simple roots of $sl_n$. 
The series (\ref{eq:kacfrenkel})
is a well-known expression of
${\rm ch}\,{\cal L}(\Lambda_k)$ \cite{kac}.
Therefore we have obtained an alternative
proof of Theorem \ref{thm:DJKMO}.

As another corollary of Theorem \ref{thm:rogers},
a new combinatorial description of the
 Kostka-Foulkes polynomials is obtained.
See Appendix \ref{sec:KF}.

\section{Vertex model of $U_q(A^{(2)}_{2n})$}

One can apply the skew diagram method also
to a vertex
model associated to the quantized
twisted affine algebra $U_q(A^{(2)}_{2n})$ \cite{kuniba}
with an appropriate modification.
We expect that 
 the underlying algebra for the degeneracy
is the Yangian of $B_n$, $Y(B_n)$.

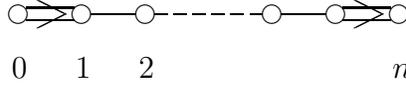
\begin{figure}[tb]
\begin{center}
\setlength{\unitlength}{1.2pt}
\begin{picture}(120,25)(0,-20)
\put(2.5,2){\line(1,0){15}}
\put(2.5,-2){\line(1,0){15}}
\put(23,0){\line(1,0){14}}
\multiput(43,0)(6,0){6}{\line(1,0){4}}
\put(83,0){\line(1,0){14}}
\put(102.5,2){\line(1,0){15}}
\put(102.5,-2){\line(1,0){15}}
\put(15,0){\line(-2,1){9}}
\put(15,0){\line(-2,-1){9}}
\put(115,0){\line(-2,1){9}}
\put(115,0){\line(-2,-1){9}}
\put(0,0){\circle{6}}
\put(20,0){\circle{6}}
\put(40,0){\circle{6}}
\put(80,0){\circle{6}}
\put(100,0){\circle{6}}
\put(120,0){\circle{6}}
\put(-2,-20){$0$}
\put(18,-20){$1$}
\put(38,-20){$2$}
\put(118,-20){$n$}
\end{picture}
\end{center}
\caption{The Dynkin diagram of $A^{(2)}_{2n}$.
The Dynkin diagram of $B_{n}$ is obtained by removing
the node $0$.}
\label{fig:dynkin}
\end{figure}

In contrast with the standard textbook \cite{kac},
we regard $A^{(2)}_{2n}$ as an affinization of
the Lie algebra $B_n$ rather than $C_n$.
Their Dynkin
diagrams are depicted in Fig.\ \ref{fig:dynkin}.
Let $\overline{\Lambda}_i$ ($i=1,\dots,n$) be
the fundamental weights of $B_n$,
and let $\epsilon_{\pm i}=\pm(\overline\Lambda_i
-\overline\Lambda_{i-1})$ for $i=1,\dots,n-1$,
$\epsilon_{\pm n}=\pm(2\overline\Lambda_n
-\overline\Lambda_{n-1})$, and $\epsilon_0=0$.
Then $B(\overline\Lambda_1)=
\{ \epsilon_1\prec\cdots\prec\epsilon_n\prec
\epsilon_0\prec\epsilon_{-n}\prec\cdots
\prec\epsilon_{-1}\}$
is the set of all the weights of the
vector representation of $B_{n}$ with
a total ordering $\prec$.
The local energy function $H:B(\overline\Lambda_1)
\times B(\overline\Lambda_1)\rightarrow
\{0,1\}$ is defined as
\begin{equation}
\label{eq:twistlocalenergy}
H(s,s')=
\begin{cases}
0 & \mbox{if $s \prec s'$ or $(s,s')=(
\epsilon_0,\epsilon_0)$}, \cr
1 & \mbox{if $s\succeq s'$ and
$(s,s')\neq (\epsilon_0,\epsilon_0)$}. \cr
\end{cases}
\end{equation}
A sequence $\vec{s}=(s_i)$, $s_i\in B(\overline\Lambda_1)$,
is a spin configuration if $\vec{s}\approx
((\epsilon_0)^\infty)$. Let ${\cal S}$ be the set of
all the spin configurations.
For each $\vec{s}\in {\cal S}$ we define 
\begin{equation}\label{eq:twistweight}
E(\vec{s})=\sum_{i=1}^\infty
i H(s_i,s_{i+1}),
\quad
{\rm wt}(\vec{s})
=-\overline\Lambda_n + \sum_{i=1}^\infty
s_i.
\end{equation}
\begin{thm}[\cite{kang}]\label{thm:kuniba}
Let ${\rm ch}\, {\cal L}(\Lambda_n)$ be
the unnormalized character of
 the (unique)
level $1$ integrable module of $A^{(2)}_{2n}$.
Then 
\begin{equation}\label{eq:atwocharacter}
{\rm ch}\, {\cal L}(\Lambda_n)(q,x)
=
\sum_{\vec{s}\in {\cal S}}
 q^{E(\vec{s})}
e^{{\rm wt}(\vec{s})},
\quad
e^{\epsilon_{\pm i}}=x_i^{\pm 1},
\ e^{\epsilon_0}=1.
\end{equation}
\noindent
\end{thm}
See Appendix \ref{sec:character} for the explicit expression
of ${\rm ch}\, {\cal L}(\Lambda_n)(q,x)$.

The local energy map $h:\vec{s}\mapsto
\vec{h}=(h_i)$, $h_i=H(s_i,s_{i+1})$,
has the image $h({\cal S})={\rm Sp}$,
where ${\rm Sp}=\{\vec{h}\mid
\vec{h}\approx ((0)^\infty)\}$.
Any element $\vec{h}\in {\rm Sp}$
is written in the form
\begin{equation*}
[m_1,\dots,m_r]:=
(\underbrace{0,\dots,0,1}_{m_1},
\dots,
\underbrace{0,\dots,0,1}_{m_r},
(0)^\infty),
\quad
m_i \geq 1.
\end{equation*}
We set ${\cal S}_{\vec{h}}=h^{-1}(\vec{h})$
for $\vec{h}\in {\rm Sp}$,
and define
$\chi_{\vec{h}}(x)
=\sum_{\vec{s}\in {\cal S}_{\vec{h}}}
e^{{\rm wt}(\vec{s})}$.

For each $\vec{h}\in {\rm Sp}$,
we associate a skew diagram
$\tilde{\kappa}(\vec{h})$
of infinite-size, following the
procedures 1--3 in the beginning
of section \ref{sect:correspondence}.
Namely, if $\vec{h}=[m_1,\dots,m_r]$,
then $\tilde\kappa(\vec{h})$
is a border strip with $r+1$ columns,
$\langle m_1,\dots,m_r,\infty
\rangle$.
Also we define a border strip
$\kappa(\vec{h})=
\langle m_1,\dots,m_r,2n\rangle$
as a ``finite part'' of 
$\tilde\kappa(\vec{h})$.

An analogue of the semi-standard condition
suited for $Y(B_n)$ 
is introduced in \cite{kunibas}
to characterize the spectrum of the
row-to-row transfer matrices
with a wide class of auxiliary spaces.
We find that the notion in \cite{kunibas} are
quite adequate also for the description of
$\chi_{\vec{h}}$.
Let $J=\{ 1 \prec \cdots \prec n
\prec 0 \prec -n \prec\cdots \prec -1\}$.
We inscribe the numbers from $J$ in
each box of a skew diagram $\lambda/\mu$.
We call such an arrangement an {\it admissible
tableau} of shape $\lambda/\mu$
if it satisfies the following condition:
Let $a$ and $b$ be the inscribed
numbers in any pair of adjacent
boxes. Then,
\begin{align}\label{eq:twistadjacent}
\mbox{(i)} &\ \mbox{
$a\prec b$ or $(a,b)=(0,0)$
 if $b$ is lower-adjacent to $a$.}
\tag{\theequation a}
\\
\mbox{(ii)} &\ \mbox{
$a\succeq b$ and $(a,b)\neq(0,0)$
if $b$ is left-adjacent to $a$.}
\tag{\theequation b}
\end{align}
Let
\addtocounter{equation}{1}
${\rm AT}(\lambda/\mu)$
be the set of all the admissible tableaux of shape 
$\lambda/\mu$.
Since the condition (\ref{eq:twistadjacent}a,b)
is just in coordinate
with (\ref{eq:twistlocalenergy}),
there is a natural bijection between 
${\cal S}_{\vec{h}}$ and 
${\rm AT}(\tilde\kappa(
\vec{h}))$ as in section \ref{sect:correspondence}.
However, the set
${\rm AT}(\kappa(
\vec{h}))$
is larger than
${\rm AT}(\tilde\kappa(
\vec{h}))$. Therefore we introduce a further
constraint on the elements in
${\rm AT}(\kappa(\vec{h}))$.
An admissible tableau of shape
$\kappa(\vec{h})$ is 
{\it L(eft)-admissible} if the content 
of the bottom $n$ boxes in the leftmost
column is frozen to the sequence
$-n, -n+1,\dots,-1$.\footnote{
The definition of the L-admissibility here
is the simplified one especially for
a border strip with the length of the
last column $2n$.
See \cite{kunibas} for the definition for a
general L-hatched skew diagram.
}
We write the set of all the
L-admissible tableaux of shape
$\kappa(\vec{h})$ as
${\rm LAT}(\kappa(\vec{h}))$.
Then
\begin{lem}
There is a bijection from ${\rm AT}
(\tilde\kappa(\vec{h}))$
to ${\rm LAT}
(\kappa(\vec{h}))$.
\end{lem}
The correspondence is a natural one:
For $\tilde{T}\in {\rm AT}
(\tilde\kappa(\vec{h}))$,
let $(a_1,\dots,a_m,(0)^\infty)$
be its content, where 
$a_{m-n+1},\dots,a_m$ are in the
top $n$ boxes in the leftmost column.
Then, $(a_1,\dots,a_m,-n,\dots,-1)$
gives the content of the corresponding
tableau
$T\in {\rm LAT}(\kappa(\vec{h}))$.
Combining the two bijections, we obtain a
bijection $\varphi:{\rm LAT}(\kappa(
\vec{h}))\rightarrow {\cal S}_{\vec{h}}$.

For an L-admissible tableau $T\in {\rm LAT}(\kappa(\vec{h}))$
with the content $(a_1, \dots, a_m, -n, \dots,-1)$,
its weight is defined as
${\rm wt}'(T)=\sum_{i=1}^m \epsilon_{a_i}+
{1\over2}\sum_{i=1}^n \epsilon_{-i}=
\sum_{i=1}^m \epsilon_{a_i}-\overline\Lambda_n$,
where, following \cite{kunibas},
we multiply the factor ${1\over2}$
on the weights corresponding to the
bottom $n$ boxes in the leftmost
 column of $\kappa(\vec{h})$.
Comparing it with (\ref{eq:twistweight}),
we see that the bijection
$\varphi:{\rm LAT}(\kappa(
\vec{h}))\rightarrow {\cal S}_{\vec{h}}$
is weight-preserving. Therefore we have
\begin{thm}
(i) The character $\chi_{\vec{h}}$ of
${\cal S}_{\vec{h}}$ is equal to
\begin{equation*}
s^{\rm L}_{\kappa(\vec{h})}:=
\sum_{T\in {\rm LAT}(\kappa(\vec{h}))}
e^{{\rm wt}'(T)}.
\end{equation*}
(ii) Let $\vec{h}=[m_1,\dots,m_r]
\in {\rm Sp}$. Then
\begin{equation*}
s^{\rm L}_{\kappa(\vec{h})}=
s^{\rm L}_{\langle m_1,\dots,m_r,2n\rangle}
=
\sigma 
\left|
\begin{array}{llllll}
1  & 1&
&\cdots&& 1\\
1 & t_{m_{r}}&t_{m_r+m_{r-1}}&&&t_{m_r+\cdots+m_1}\\
0 & 1&t_{m_{r-1}} \\
&0&1&&&\vdots\\
&&&\ddots\\
&&&1& t_{m_2}&t_{m_2+m_1}\\ 
&&&0&1& t_{m_1} 
\end{array}
\right|\,,
\end{equation*}
where $\sigma$ is ${\rm ch}\, 
V_{\overline\Lambda_n}$,
$t_m$ is $0$ for\/ $m<0$,
${\rm ch}\, (V_{\overline\Lambda_m}
\oplus V_{\overline\Lambda_{m-2}}\oplus
\cdots)$ for\/ $0\leq m\leq n-1$, and
$\sigma^2-t_{2n-1-m}$ for\/ $m\geq n$,
and $V_{\overline\Lambda_j}$ is the $j${\rm th}
fundamental representation of $B_n$.
\hfill\break
(iii) The character of the level $1$
integrable module ${\cal L}(\Lambda_n)$
of $A^{(2)}_{2n}$ decomposes as
\begin{equation}\label{eq:twistch}
{\rm ch}\, {\cal L}(\Lambda_n)(q,x)
=
\sum_{\kappa\in  BS}
q^{t(\kappa)} 
s^{\rm L}_{\kappa}(x),
\end{equation}
where BS is the set of all the border strips
$\kappa=\langle m_1,\dots,m_r\rangle$
with $m_r=2n$, and $t(\kappa)=
\sum_{i=1}^{r-1}(r-i)m_i$.
\end{thm}
\begin{pf}
(i) and (iii) are due to the bijection
$\varphi$ and Theorem \ref{thm:kuniba}.
(ii) is a special case of Theorem 4.1
in \cite{kunibas}, which is an analogue
of Prop.\ \ref{prop:equivalent} for
$s^{\rm L}_{\lambda/\mu}$.
\end{pf}

It was conjectured in \cite{kunibas}
that $s^{\rm L}_{\kappa(\vec{h})}$ is
the $B_n$-character of a certain
irreducible $Y(B_n)$-module.
Therefore it is natural to conjecture that
\begin{conj}
There is a canonical action of $Y(B_n)$ on
the level $1$ integrable module 
${\cal L}(\Lambda_n)$ of
$A^{(2)}_{2n}$.
The decomposition (\ref{eq:twistch})
describes its $Y(B_n)$-module structure
on ${\cal L}(\Lambda_n)$.
\end{conj}

\section{Conclusion}

In this paper we exhibit intimate
relationships among the 
spectral decomposition of the
vertex models, skew diagrams and
the associated Schur functions,
irreducible characters of the Yangians,
Yangian module structures in conformal
field theory,
spectra of spin models with the
inverse-square interaction, and so on.
We believe that further study of this 
interrelation will enlighten our
understanding of the common integrable
structure behind these models.

It is also interesting to investigate
other vertex models.
For example,
for the symmetric fusion models of
$U_q(\widehat{sl}_n)$,
which correspond to the higher level 
integrable modules of $\widehat{sl}_n$,
it is possible to extend our skew diagram
approach. 
Even though conceptually it is
quite analogous to the level 1 case,
some 
complexity enters, especially for $n\geq3$.
Notable changes are,
firstly, skew diagrams of non-border
strips are necessary to describe the spectrum,
and secondly, the characters of
non-tame modules appear
as the characters of the degeneracy of the
spectrum. We hope to give a full report
on it in a future publication.

At the very last stage
of  the preparation of the manuscript,
the preprint ``The $\widehat{SU}(n)_1$ WZW models,
spinon decomposition and Yangian structure''
by P.\ Bouwknegt and K.\ Schoutens (hep-th/9607064)
appeared, where the authors obtain a partially
overlapping result to ours.

\vskip\baselineskip

{\it Acknowledgment.}
A.\ N.\ K.\ would like to thank the colleagues
from Tokyo University for kind hospitality
and support.
A.\ K.\ and T.\ N.\ appreciate
the great hospitality of the organizers
of
the {\it Third International Conference
on Conformal Filed Theory and Integrable Models\/}
held at Landau Institute of Theoretical Physics,
Chernogolovka on June 24--29, 1996,
where the main result of this paper is presented.
T.\ N.\ would like to thank G.\ Felder and
A.\ Varchenko for their warm hospitality
at University of North Carolina.

\appendix
\section{A proof of Theorem \ref{thm:rogers}}

In this appendix we give a proof of Theorem
\ref{thm:rogers}.

Any element of $\vec{h}\in {\rm Sp}_N$
is uniquely written as
\begin{equation*}
\vec{h}=
(\underbrace{0,\dots,0,1}_{m_1},
\underbrace{0,\dots,0,1}_{m_2},
\dots,
\underbrace{0,\dots,0,1}_{m_r},
(\underbrace{0,\dots,0,1}_{n})^\infty
)
\end{equation*}
for an integer $r$ (= the number of
1's in the first $N$ elements of $\vec{h}$)
and the integers $1\leq m_i\leq n$ such that
$m_1+m_2+\cdots +m_r=N$.
For such an  $\vec{h}\in {\cal S}_N$
we associate a border strip
$\langle m_1,\dots,m_r\rangle$.
If $m_r\neq n$,
the border strip
 $\langle m_1,\dots,m_r\rangle$ is 
equal to $\kappa (\vec{h})$.
If $m_r=n$, however,
 $\langle m_1,\dots,m_r\rangle$ has a few more
length-$n$ columns than $\kappa(\vec{h})$
on its tail such that the size of 
$\langle m_1,\dots,m_r\rangle$ is always $N$.
In either case,
we have
$s_{\kappa(\vec{h})}
=s_{\langle m_1,\dots,m_r\rangle}$
thanks to the specialization $x_1\cdots x_n=1$.
Furthermore, for $\vec{h}\in {\rm Sp}_N$
it holds that
\begin{equation*}
\sum_{i=1}^\infty
i(h_i-h_i^{(\overline{N})})=
\sum_{i=1}^{N-1}i(h_i-1) + E_N
=
\sum_{i=1}^{N}i(h_i-1) + E_N.
\end{equation*}
Thus Theorem \ref{thm:rogers} is equivalent to
\begin{thm}\label{thm:schurrogers}
\begin{equation}\label{eq:schurrogers}
\sum_{r=1}^N
\sum_{1\leq m_i \leq n\atop
m_1+\cdots + m_r=N}
q^{{1\over2}N(N+1)-
\sum_{i=1}^r (m_1+\cdots+m_i)}
s_{\langle m_1,\dots,m_r\rangle}(x)
=
 H_N(q,x).
\end{equation}
\end{thm}

\begin{rem}
Theorem \ref{thm:schurrogers}
is true without assuming the relation
$x_1\cdots x_n=1$ as we show below.
\end{rem}

We write the left hand side
 of (\ref{eq:schurrogers})
as $F_N(q,x)$ 
($F_0(q,x)=1$ by definition).
\begin{exm}
\begin{alignat*}{2}
F_0(q,x)&=1, &\qquad H_0(q,x)&=1,\\
F_1(q,x)&=s\,_{
\setlength{\unitlength}{1.2pt}
\begin{picture}(4,4)(0,0)
\put(0,0){\line(1,0){4}}
\put(0,4){\line(1,0){4}}
\put(0,0){\line(0,1){4}}
\put(4,0){\line(0,1){4}}
\end{picture}
}\,(x),
&\qquad
H_1(q,x)&=\sum_{i=1}^n x_i,\\
F_2(q,x)&=q 
s\,_{
\setlength{\unitlength}{1.2pt}
\begin{picture}(4,8)(0,0)
\put(0,0){\line(1,0){4}}
\put(0,4){\line(1,0){4}}
\put(0,8){\line(1,0){4}}
\put(0,0){\line(0,1){8}}
\put(4,0){\line(0,1){8}}
\end{picture}
}\, (x)
+
s\,_{
\setlength{\unitlength}{1.2pt}
\begin{picture}(8,4)(0,0)
\put(0,0){\line(1,0){8}}
\put(0,4){\line(1,0){8}}
\put(0,0){\line(0,1){4}}
\put(4,0){\line(0,1){4}}
\put(8,0){\line(0,1){4}}
\end{picture}
}\,(x),
&\qquad
H_2(q,x)&=\sum_{i=1}^n x_i^2
+(1+q)\sum_{1\leq i<j\leq n} x_ix_j.
\end{alignat*}
\end{exm}

Following \cite{hikami},
 we consider the recursion relation
for $F_N(q,x)$ and $H_N(q,x)$.
Let
\begin{equation*}
G(q,x,t)={1\over (tx_1;q)_\infty
(tx_2;q)_\infty \cdots
(tx_n;q)_\infty},
\quad
(a;q)_\infty = \prod_{j=0}^\infty
(1-aq^j).
\end{equation*}

\begin{lem}[\cite{andrews}]\label{lemma:generate}
The function $G(q,x,t)$ is the generating
function of $H_N(q,x)$:
\begin{equation*}
G(q,x,t)=
\sum_{N=0}^\infty
{H_N(q,x)\over (q)_N} t^N.
\end{equation*}
\end{lem}
\begin{pf}
It easily follows from the identity
\cite{andrews}
\begin{equation*}
{1\over (t;q)_\infty}
=
\sum_{N=0}^\infty {1\over (q)_N} t^N.
\end{equation*}
\end{pf}

\begin{lem}[\cite{andrews,hikami}]
The functions $H_N(q,x)$ satisfy the
following recursion relation: For any $N\geq 1$,
\begin{equation}\label{eq:recursion}
H_N(q,x)=\sum_{i=1}^n
(-1)^{i+1} {(q)_{N-1}\over
(q)_{N-i}}
e_i(x) H_{N-i}(q,x)=0,
\end{equation}
where $e_i(x)$ is the $i${\rm th} elementary symmetric
function of $x_1,\dots,x_n$,
and  $H_N(q,x)=0$ for $N<0$.
\end{lem}
\begin{pf}
Consider the identity,
\begin{equation*}
G(q,x,qt)=
\left(\prod_{i=1}^n(1-tx_i)
\right)
G(q,x,t)
=
\left(\sum_{i=0}^n
(-1)^i e_i(x) t^i\right)
G(q,x,t).
\end{equation*}
Comparing the coefficients of $t^N$ of
the both sides, and using Lemma
 \ref{lemma:generate}, we have
\begin{equation*}
{q^N\over (q)_N}
H_N(q,x)
=
\sum_{i=0}^n(-1)^i e_i(x)
{H_{N-i}(q,x)\over
(q)_{N-i}},
\end{equation*}
from which the lemma follows.
\end{pf}

The recursion relation (\ref{eq:recursion}),
together
with the initial condition $H_0(q,x)=1$,
uniquely determines $H_N(q,x)$.
In the rest of the appendix  we show 
$F_N(q,x)$ also satisfies (\ref{eq:recursion}).
We recall
\begin{equation*}
s_{\langle m_1,\dots,m_r\rangle}
=
\left|
\begin{array}{llllll}
e_{m_r} & e_{m_r+m_{r-1}}&
&\cdots&& e_{m_r+\cdots+m_1}\\
1 & e_{m_{r-1}}\\
0 & 1\\
&0&&&&\vdots\\
&&\ddots&\ddots\\
&&0&1& e_{m_2}&e_{m_2+m_1}\\ 
&&&0&1& e_{m_1} 
\end{array}
\right|\,.
\end{equation*}
Expanding the determinant along the first
row, we have
\begin{lem}
\begin{equation}\label{eq:schurrec}
s_{\langle m_1,\dots,m_r\rangle}=
\sum_{i=1}^r(-1)^{i+1}
e_{m_r+\cdots + m_{r-i+1}}
s_{\langle m_1,\dots,m_{r-i}\rangle}.
\end{equation}
\end{lem}
Substituting (\ref{eq:schurrec})
into $F_N(q,x)$, we have
\begin{lem}
The functions $F_N(q,x)$ satisfy
the following recursion relation: For any $N\geq 1$,
\begin{align}\label{eq:frecursion}
F_N(q,x)
&=\sum_{m=1}^N A_{N,m}(q)
e_m(x) F_{N-m}(q,x),\\
A_{N,m}(q)
&=\sum_{j=1}^m
\sum_{1\leq k_i\leq n\atop
k_1+\cdots + k_j=m}
(-1)^{j+1} q^{c_{N,m}(k_1,\dots,k_j)},\\
c_{N,m}(k_1,\dots,k_j)
&=
{1\over2}N(N+1)-{1\over2}(N-m)(N-m+1)
\notag\\
&\quad
-\sum_{i=1}^j (N-m+k_1+\cdots + k_i)
\notag\\
&=
N(m-j)-{1\over2}m^2-{1\over2}m+
\sum_{i=1}^j i k_i,
\end{align}
where $F_N(q,x)=0$ for $N<0$.
\end{lem}

In the right hand side of
(\ref{eq:frecursion}) we can replace
the summation $\sum_{m=1}^N$ by
$\sum_{m=1}^n$ because $e_m(x)=0$
for $m>n$.
To complete the proof of
Theorem \ref{thm:schurrogers}, we have
only to show that

\begin{lem}
For $1\leq m \leq {\rm min}(n,N)$,
\begin{equation}\label{eq:afunction}
A_{N,m}(q)=(-1)^{m+1}
{(q)_{N-1}\over(q)_{N-m}}.
\end{equation}
\end{lem}
\begin{pf}
For any $N\geq 1$,
$A_{N,1}(q)=q^{c_{N,1}(1)}=1$.
Below we show
\begin{equation}\label{eq:arecursion}
A_{N,m}(q)=-(1-q^{N-1})A_{N-1,m-1}(q),
\end{equation}
for $2\leq m \leq N$. 
{}From these, the lemma follows.

Consider
\begin{equation}\label{eq:asum}
A_{N,m}(q)=
\sum_{j=1}^{m}
\sum_{k_i\geq 1\atop
k_1+\cdots + k_j=m}
(-1)^{j+1}
q^{c_{N,m}(k_1,\dots,k_j)}
\end{equation}
for $2\leq m\leq N$. Notice that we have dropped the upper
inequality $k_i \leq n$ in the summation, because it
is automatically satisfied under the assumption
 $m\leq n$.
Let 
\begin{equation*}
I_m=\{ (k_1,\dots,k_j)\mid
j,k_i\geq 1, k_1+\cdots +k_j=m\}
\end{equation*}
be the set of all the ordered partitions
of $m$.
Then $I_{m}$ is the disjoint union
of $I_{m}^{(1)}$ and $I_{m}^{(2)}$,
where
\begin{align*}
I_{m}^{(1)}
&= \{ (k_1,\dots,k_j,1)\mid
(k_1,\dots,k_j)\in I_{m-1}\},\\
I_{m}^{(2)}
&= \{ (k_1,\dots,k_j+1)\mid
(k_1,\dots,k_j)\in I_{m-1}\}.
\end{align*}

Let us perform the summation in (\ref{eq:asum})
over $I_{m}^{(1)}$ and $I_{m}^{(2)}$,
separately.
The contribution from $I_{m}^{(1)}$
is
\begin{align*}
&\sum_{j=1}^{m-1}
\sum_{(k_1,\dots,k_j)\in I_{m-1}}
(-1)^{j+2} q^{c_{N,m}(k_1,\dots,k_j,1)}\\
=&-\sum_{j=1}^{m-1}
\sum_{(k_1,\dots,k_j)\in I_{m-1}}
(-1)^{j+1} q^{c_{N-1,m-1}(k_1,\dots,k_j)}
=-A_{N-1,m-1}(q).
\end{align*}
The contribution from $I_{m}^{(2)}$
is
\begin{align*}
&\sum_{j=1}^{m-1}
\sum_{(k_1,\dots,k_j)\in I_{m-1}}
(-1)^{j+1} q^{c_{N,m}(k_1,\dots,k_j+1)}\\
=&\sum_{j=1}^{m-1}
\sum_{(k_1,\dots,k_j)\in I_{m-1}}
(-1)^{j+1} q^{c_{N-1,m-1}(k_1,\dots,k_j)+N-1}
=
q^{N-1}A_{N-1,m-1}(q).
\end{align*}
Putting the both together, we get
(\ref{eq:arecursion}).
\end{pf}

\section{A new combinatorial formula for Kostka-Foulkes
polynomials}\label{sec:KF}

As another corollary of Theorem \ref{thm:rogers},
or Theorem
\ref{thm:schurrogers},  we get a new formula for 
the Kostka-Foulkes polynomials in terms of the Littlewood-Richardson
tableaux of border strips.
See \cite{kirillov,macdonald} for further information
of the material discussed here.
For a border strip  $\kappa=\langle m_1,\dots,m_r\rangle$
we define  $t(\kappa)$ as
\begin{equation*}
t(\kappa) = \sum_{i=1}^{r-1} (r-i)m_i.
\end{equation*}
Further,  for a partition $\lambda$,
we denote by  $C(\kappa,\lambda)$
the number of the semi-standard tableaux that form 
lattice permutations (cf.\ \cite{macdonald})
and are of shape $\kappa$ and content $\lambda$.
The number $C(\kappa=\mu/\nu,\lambda)$ has an interpretation
as $c^\mu_{\nu\lambda}=
{\rm Mult}_{V_{\lambda}}(V_{\mu,\nu})$, i.e., 
the multiplicity of the irreducible 
representation $V_\lambda$ of $sl_n$ in the restriction of 
the representation $V_{\mu,\nu}$ of $Y(sl_n)$ to $sl_n$.
\begin{prop}
Let $\lambda$ be a partition, then 
\begin{equation}\label{eq:kostka}
\sum_\kappa q^{t(\kappa)} C(\kappa, \lambda)
 = K_{\lambda, (1^{\vert \lambda \vert})}(q),
\end{equation}
where the summation is taken over all the border strips  $\kappa$
of rank $n$ with $\vert \kappa \vert = \vert \lambda \vert$.
\end{prop}
\begin{pf}
Theorem \ref{thm:schurrogers} can be rewritten as
\begin{equation*}
\sum_{\kappa:{\rm border\ strip}\atop |\kappa|=N}
q^{t(\kappa)}s_\kappa(x)
=
\sum_{k_i\in {\bf Z}_{\geq 0}\atop
k_1+\cdots +k_n=N}
q^{\sum_{i=1}^n {1\over2} k_i (k_i-1)}
{(q)_N\over (q)_{k_1}\cdots (q)_{k_n}}
x_1^{k_1}\cdots x_n^{k_n}.
\end{equation*}
On the other hand, the degree-$N$ part of Corollary 6
of \cite{kirillov} reads as
\begin{equation*}
\sum_{\lambda:{\rm Young\ diagram}\atop |\lambda|=N}
K_{\lambda,(1^{|\lambda|})}(q)s_\lambda(x)
=
\sum_{k_i\in {\bf Z}_{\geq 0}\atop
k_1+\cdots +k_n=N}
q^{\sum_{i=1}^n {1\over2} k_i (k_i-1)}
{(q)_N\over (q)_{k_1}\cdots (q)_{k_n}}
x_1^{k_1}\cdots x_n^{k_n}.
\end{equation*}
Equate the left hand sides of two equalities,
and expand $s_{\kappa}(x)$ as
$\sum_\lambda C(\kappa,\lambda)
 s_\lambda(x)$, where 
$|\lambda|=|\kappa|$ in the summation.
Using the independence of $s_\lambda$'s,
 we have (\ref{eq:kostka}).
\end{pf}

\begin{exm}
Let $\lambda = (3,2,1)$. Then it is known \cite{macdonald}
 that
\begin{equation*}
K_{(3,2,1), (1^6)}(q) = q^4(1+q)^2(1+q^2)(1+q^3).
\end{equation*}
On the other hand, there exist  14 border strips
which give non-zero contribution to the left
hand side of (\ref{eq:kostka}) with the following tableaux:
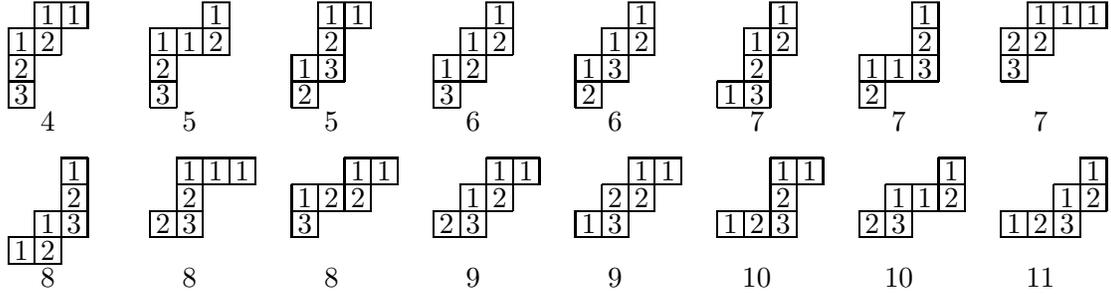
\small
\begin{center}
\begin{tabular}{rrrrrrrr}
\setlength{\unitlength}{1pt}
\begin{picture}(40,50)(0,-50)
\put(0,-40){\line(1,0){10}}
\put(0,-30){\line(1,0){10}}
\put(0,-20){\line(1,0){20}}
\put(0,-10){\line(1,0){30}}
\put(10,0){\line(1,0){20}}
\put(30,0){\line(0,-1){10}}
\put(20,0){\line(0,-1){20}}
\put(10,0){\line(0,-1){40}}
\put(0,-10){\line(0,-1){30}}
\put(20,-10){{\makebox(10,10){$1$}}}
\put(10,-10){{\makebox(10,10){$1$}}}
\put(10,-20){{\makebox(10,10){$2$}}}
\put(0,-20){{\makebox(10,10){$1$}}}
\put(0,-30){{\makebox(10,10){$2$}}}
\put(0,-40){{\makebox(10,10){$3$}}}
\put(10,-50){{\makebox(10,10){$4$}}}
\end{picture}
&
\setlength{\unitlength}{1pt}
\begin{picture}(40,50)(0,-50)
\put(0,-40){\line(1,0){10}}
\put(0,-30){\line(1,0){10}}
\put(0,-20){\line(1,0){30}}
\put(0,-10){\line(1,0){30}}
\put(20,0){\line(1,0){10}}
\put(30,0){\line(0,-1){20}}
\put(20,0){\line(0,-1){20}}
\put(10,-10){\line(0,-1){30}}
\put(0,-10){\line(0,-1){30}}
\put(20,-10){{\makebox(10,10){$1$}}}
\put(10,-20){{\makebox(10,10){$1$}}}
\put(20,-20){{\makebox(10,10){$2$}}}
\put(0,-20){{\makebox(10,10){$1$}}}
\put(0,-30){{\makebox(10,10){$2$}}}
\put(0,-40){{\makebox(10,10){$3$}}}
\put(10,-50){{\makebox(10,10){$5$}}}
\end{picture}
&
\setlength{\unitlength}{1pt}
\begin{picture}(40,50)(0,-50)
\put(0,-40){\line(1,0){10}}
\put(0,-30){\line(1,0){20}}
\put(0,-20){\line(1,0){20}}
\put(10,-10){\line(1,0){20}}
\put(10,0){\line(1,0){20}}
\put(30,0){\line(0,-1){10}}
\put(20,0){\line(0,-1){30}}
\put(10,0){\line(0,-1){40}}
\put(0,-20){\line(0,-1){20}}
\put(20,-10){{\makebox(10,10){$1$}}}
\put(10,-10){{\makebox(10,10){$1$}}}
\put(0,-30){{\makebox(10,10){$1$}}}
\put(10,-20){{\makebox(10,10){$2$}}}
\put(0,-40){{\makebox(10,10){$2$}}}
\put(10,-30){{\makebox(10,10){$3$}}}
\put(10,-50){{\makebox(10,10){$5$}}}
\end{picture}
&
\setlength{\unitlength}{1pt}
\begin{picture}(40,50)(0,-50)
\put(0,-40){\line(1,0){10}}
\put(0,-30){\line(1,0){20}}
\put(0,-20){\line(1,0){30}}
\put(10,-10){\line(1,0){20}}
\put(20,0){\line(1,0){10}}
\put(30,0){\line(0,-1){20}}
\put(20,0){\line(0,-1){30}}
\put(10,-10){\line(0,-1){30}}
\put(0,-20){\line(0,-1){20}}
\put(20,-10){{\makebox(10,10){$1$}}}
\put(10,-20){{\makebox(10,10){$1$}}}
\put(0,-30){{\makebox(10,10){$1$}}}
\put(20,-20){{\makebox(10,10){$2$}}}
\put(10,-30){{\makebox(10,10){$2$}}}
\put(0,-40){{\makebox(10,10){$3$}}}
\put(10,-50){{\makebox(10,10){$6$}}}
\end{picture}
&
\setlength{\unitlength}{1pt}
\begin{picture}(40,50)(0,-50)
\put(0,-40){\line(1,0){10}}
\put(0,-30){\line(1,0){20}}
\put(0,-20){\line(1,0){30}}
\put(10,-10){\line(1,0){20}}
\put(20,0){\line(1,0){10}}
\put(30,0){\line(0,-1){20}}
\put(20,0){\line(0,-1){30}}
\put(10,-10){\line(0,-1){30}}
\put(0,-20){\line(0,-1){20}}
\put(20,-10){{\makebox(10,10){$1$}}}
\put(10,-20){{\makebox(10,10){$1$}}}
\put(0,-30){{\makebox(10,10){$1$}}}
\put(20,-20){{\makebox(10,10){$2$}}}
\put(0,-40){{\makebox(10,10){$2$}}}
\put(10,-30){{\makebox(10,10){$3$}}}
\put(10,-50){{\makebox(10,10){$6$}}}
\end{picture}
&
\setlength{\unitlength}{1pt}
\begin{picture}(40,50)(0,-50)
\put(0,-40){\line(1,0){20}}
\put(0,-30){\line(1,0){20}}
\put(10,-20){\line(1,0){20}}
\put(10,-10){\line(1,0){20}}
\put(20,0){\line(1,0){10}}
\put(30,0){\line(0,-1){20}}
\put(20,0){\line(0,-1){40}}
\put(10,-10){\line(0,-1){30}}
\put(0,-30){\line(0,-1){10}}
\put(20,-10){{\makebox(10,10){$1$}}}
\put(10,-20){{\makebox(10,10){$1$}}}
\put(0,-40){{\makebox(10,10){$1$}}}
\put(20,-20){{\makebox(10,10){$2$}}}
\put(10,-30){{\makebox(10,10){$2$}}}
\put(10,-40){{\makebox(10,10){$3$}}}
\put(10,-50){{\makebox(10,10){$7$}}}
\end{picture}
&
\setlength{\unitlength}{1pt}
\begin{picture}(40,50)(0,-50)
\put(0,-40){\line(1,0){10}}
\put(0,-30){\line(1,0){30}}
\put(0,-20){\line(1,0){30}}
\put(20,-10){\line(1,0){10}}
\put(20,0){\line(1,0){10}}
\put(30,0){\line(0,-1){30}}
\put(20,0){\line(0,-1){30}}
\put(10,-20){\line(0,-1){20}}
\put(0,-20){\line(0,-1){20}}
\put(20,-10){{\makebox(10,10){$1$}}}
\put(10,-30){{\makebox(10,10){$1$}}}
\put(0,-30){{\makebox(10,10){$1$}}}
\put(20,-20){{\makebox(10,10){$2$}}}
\put(0,-40){{\makebox(10,10){$2$}}}
\put(20,-30){{\makebox(10,10){$3$}}}
\put(10,-50){{\makebox(10,10){$7$}}}
\end{picture}
&
\setlength{\unitlength}{1pt}
\begin{picture}(40,50)(0,-50)
\put(0,-30){\line(1,0){10}}
\put(0,-20){\line(1,0){20}}
\put(0,-10){\line(1,0){40}}
\put(10,0){\line(1,0){30}}
\put(40,0){\line(0,-1){10}}
\put(30,0){\line(0,-1){10}}
\put(20,0){\line(0,-1){20}}
\put(10,0){\line(0,-1){30}}
\put(0,-10){\line(0,-1){20}}
\put(30,-10){{\makebox(10,10){$1$}}}
\put(20,-10){{\makebox(10,10){$1$}}}
\put(10,-10){{\makebox(10,10){$1$}}}
\put(10,-20){{\makebox(10,10){$2$}}}
\put(0,-20){{\makebox(10,10){$2$}}}
\put(0,-30){{\makebox(10,10){$3$}}}
\put(10,-50){{\makebox(10,10){$7$}}}
\end{picture}
\end{tabular}
\end{center}


\begin{center}
\begin{tabular}{rrrrrrrr}
\setlength{\unitlength}{1pt}
\begin{picture}(40,50)(0,-50)
\put(0,-40){\line(1,0){20}}
\put(0,-30){\line(1,0){30}}
\put(10,-20){\line(1,0){20}}
\put(20,-10){\line(1,0){10}}
\put(20,0){\line(1,0){10}}
\put(30,0){\line(0,-1){30}}
\put(20,0){\line(0,-1){40}}
\put(10,-20){\line(0,-1){20}}
\put(0,-30){\line(0,-1){10}}
\put(20,-10){{\makebox(10,10){$1$}}}
\put(10,-30){{\makebox(10,10){$1$}}}
\put(0,-40){{\makebox(10,10){$1$}}}
\put(20,-20){{\makebox(10,10){$2$}}}
\put(10,-40){{\makebox(10,10){$2$}}}
\put(20,-30){{\makebox(10,10){$3$}}}
\put(10,-50){{\makebox(10,10){$8$}}}
\end{picture}
&
\setlength{\unitlength}{1pt}
\begin{picture}(40,50)(0,-50)
\put(0,-30){\line(1,0){20}}
\put(0,-20){\line(1,0){20}}
\put(10,-10){\line(1,0){30}}
\put(10,0){\line(1,0){30}}
\put(40,0){\line(0,-1){10}}
\put(30,0){\line(0,-1){10}}
\put(20,0){\line(0,-1){30}}
\put(10,0){\line(0,-1){30}}
\put(0,-20){\line(0,-1){10}}
\put(30,-10){{\makebox(10,10){$1$}}}
\put(20,-10){{\makebox(10,10){$1$}}}
\put(10,-10){{\makebox(10,10){$1$}}}
\put(10,-20){{\makebox(10,10){$2$}}}
\put(0,-30){{\makebox(10,10){$2$}}}
\put(10,-30){{\makebox(10,10){$3$}}}
\put(10,-50){{\makebox(10,10){$8$}}}
\end{picture}
&
\setlength{\unitlength}{1pt}
\begin{picture}(40,50)(0,-50)
\put(0,-30){\line(1,0){10}}
\put(0,-20){\line(1,0){30}}
\put(0,-10){\line(1,0){40}}
\put(20,0){\line(1,0){20}}
\put(40,0){\line(0,-1){10}}
\put(30,0){\line(0,-1){20}}
\put(20,0){\line(0,-1){20}}
\put(10,-10){\line(0,-1){20}}
\put(0,-10){\line(0,-1){20}}
\put(30,-10){{\makebox(10,10){$1$}}}
\put(20,-10){{\makebox(10,10){$1$}}}
\put(0,-20){{\makebox(10,10){$1$}}}
\put(20,-20){{\makebox(10,10){$2$}}}
\put(10,-20){{\makebox(10,10){$2$}}}
\put(0,-30){{\makebox(10,10){$3$}}}
\put(10,-50){{\makebox(10,10){$8$}}}
\end{picture}
&
\setlength{\unitlength}{1pt}
\begin{picture}(40,50)(0,-50)
\put(0,-30){\line(1,0){20}}
\put(0,-20){\line(1,0){30}}
\put(10,-10){\line(1,0){30}}
\put(20,0){\line(1,0){20}}
\put(40,0){\line(0,-1){10}}
\put(30,0){\line(0,-1){20}}
\put(20,0){\line(0,-1){30}}
\put(10,-10){\line(0,-1){20}}
\put(0,-20){\line(0,-1){10}}
\put(30,-10){{\makebox(10,10){$1$}}}
\put(20,-10){{\makebox(10,10){$1$}}}
\put(10,-20){{\makebox(10,10){$1$}}}
\put(20,-20){{\makebox(10,10){$2$}}}
\put(0,-30){{\makebox(10,10){$2$}}}
\put(10,-30){{\makebox(10,10){$3$}}}
\put(10,-50){{\makebox(10,10){$9$}}}
\end{picture}
&
\setlength{\unitlength}{1pt}
\begin{picture}(40,50)(0,-50)
\put(0,-30){\line(1,0){20}}
\put(0,-20){\line(1,0){30}}
\put(10,-10){\line(1,0){30}}
\put(20,0){\line(1,0){20}}
\put(40,0){\line(0,-1){10}}
\put(30,0){\line(0,-1){20}}
\put(20,0){\line(0,-1){30}}
\put(10,-10){\line(0,-1){20}}
\put(0,-20){\line(0,-1){10}}
\put(30,-10){{\makebox(10,10){$1$}}}
\put(20,-10){{\makebox(10,10){$1$}}}
\put(0,-30){{\makebox(10,10){$1$}}}
\put(20,-20){{\makebox(10,10){$2$}}}
\put(10,-20){{\makebox(10,10){$2$}}}
\put(10,-30){{\makebox(10,10){$3$}}}
\put(10,-50){{\makebox(10,10){$9$}}}
\end{picture}
&
\setlength{\unitlength}{1pt}
\begin{picture}(40,50)(0,-50)
\put(0,-30){\line(1,0){30}}
\put(0,-20){\line(1,0){30}}
\put(20,-10){\line(1,0){20}}
\put(20,0){\line(1,0){20}}
\put(40,0){\line(0,-1){10}}
\put(30,0){\line(0,-1){30}}
\put(20,0){\line(0,-1){30}}
\put(10,-20){\line(0,-1){10}}
\put(0,-20){\line(0,-1){10}}
\put(30,-10){{\makebox(10,10){$1$}}}
\put(20,-10){{\makebox(10,10){$1$}}}
\put(0,-30){{\makebox(10,10){$1$}}}
\put(20,-20){{\makebox(10,10){$2$}}}
\put(10,-30){{\makebox(10,10){$2$}}}
\put(20,-30){{\makebox(10,10){$3$}}}
\put(10,-50){{\makebox(10,10){$10$}}}
\end{picture}
&
\setlength{\unitlength}{1pt}
\begin{picture}(40,50)(0,-50)
\put(0,-30){\line(1,0){20}}
\put(0,-20){\line(1,0){40}}
\put(10,-10){\line(1,0){30}}
\put(30,0){\line(1,0){10}}
\put(40,0){\line(0,-1){20}}
\put(30,0){\line(0,-1){20}}
\put(20,-10){\line(0,-1){20}}
\put(10,-10){\line(0,-1){20}}
\put(0,-20){\line(0,-1){10}}
\put(30,-10){{\makebox(10,10){$1$}}}
\put(20,-20){{\makebox(10,10){$1$}}}
\put(10,-20){{\makebox(10,10){$1$}}}
\put(30,-20){{\makebox(10,10){$2$}}}
\put(0,-30){{\makebox(10,10){$2$}}}
\put(10,-30){{\makebox(10,10){$3$}}}
\put(10,-50){{\makebox(10,10){$10$}}}
\end{picture}
&
\setlength{\unitlength}{1pt}
\begin{picture}(40,50)(0,-50)
\put(0,-30){\line(1,0){30}}
\put(0,-20){\line(1,0){40}}
\put(20,-10){\line(1,0){20}}
\put(30,0){\line(1,0){10}}
\put(40,0){\line(0,-1){20}}
\put(30,0){\line(0,-1){30}}
\put(20,-10){\line(0,-1){20}}
\put(10,-20){\line(0,-1){10}}
\put(0,-20){\line(0,-1){10}}
\put(30,-10){{\makebox(10,10){$1$}}}
\put(20,-20){{\makebox(10,10){$1$}}}
\put(0,-30){{\makebox(10,10){$1$}}}
\put(30,-20){{\makebox(10,10){$2$}}}
\put(10,-30){{\makebox(10,10){$2$}}}
\put(20,-30){{\makebox(10,10){$3$}}}
\put(10,-50){{\makebox(10,10){$11$}}}
\end{picture}
\end{tabular}
\end{center}

\normalsize
\noindent
where the number $t(\kappa)$ is attached below
each diagram.
Hence the left hand side of (\ref{eq:kostka})
is
\begin{equation*}
q^4(1+2q+
2q^2 + 3q^3+
3q^4+2q^5
+2q^6+q^7)=
K_{(3,2,1), (1^6)}(q).
\end{equation*}
\end{exm}

In the same spirit we define the branching function
$b^{\Lambda_k}_\lambda(q)$ of ${\cal L}(\Lambda_k)$
(as an $sl_n$-module) as
\begin{equation*}
{\rm ch}\, {\cal L}(\Lambda_k)(q,x)
=\sum_{\lambda:{\rm Young\, diagram}
\atop l(\lambda)<n,\, |\lambda|\equiv k\,
{\rm mod}\, n}
b^{\Lambda_k}_\lambda(q)
s_{\lambda}(x).
\end{equation*}
Then from (\ref{eq:characterdecomposition}a)
we have
\begin{prop}
\begin{equation*}
b^{\Lambda_k}_\lambda(q)=
q^{-{c\over24}}
\sum_{\kappa\in BS\atop
|\kappa|\geq |\lambda|,\, |\kappa|\equiv k\, {\rm mod}\, n}
 q^{{1\over2n}|\kappa|(n-|\kappa|)+t(\kappa)}
C\left(\kappa,\lambda+\left(\left({|\kappa|-|\lambda|\over n}
\right)^n\right)\right).
\end{equation*}
\end{prop}

We are going to prove an analogue of (\ref{eq:kostka})
 for the 
Kostka-Foulkes polynomials $K_{\lambda, (\ell^N)}(q)$
(which should correspond to the $Y(sl_n)$-module
structure on ${\cal L}(\ell \Lambda_0)$) in a separate publication.

\section{The character of the level $1$ integrable module of 
$A^{(2)}_{2n}$}\label{sec:character}

In this appendix we describe the function
${\rm ch}\, {\cal L}(\Lambda_n)(q,x)$ in
 Theorem \ref{thm:kuniba}.\footnote{It is our pleasure to
thank M.\ Wakimoto for his helpful comments.}

Let $\alpha_0,\dots,\alpha_n$ be the simple roots of
$A^{(2)}_{2n}$ with the label in Fig.\ \ref{fig:dynkin}.
We have the {\it null root} $\delta=\alpha_0+
2(\alpha_1+\cdots \alpha_n)$. Let $\{\epsilon_1,\dots,
\epsilon_n\}$ be the orthonormal basis of the
dual of the Cartan subalgebra of $B_n$ with
\begin{equation}\label{eq:broot}
\alpha_0=\delta-2\epsilon_1, \quad
\alpha_1=\epsilon_1-\epsilon_2,\quad\dots,\quad
\alpha_{n-1}=\epsilon_{n-1}-\epsilon_n,\quad
\alpha_n=\epsilon_n.
\end{equation}
On the other hand, let $\alpha'_i=\alpha_{n-i}$,
which favors the subalgebra $C_n$.
(Below the prime symbol $'$ indicates that we are in the $C_n$ picture.)
 Let $\{\epsilon'_1,\dots,
\epsilon'_n\}$ be the orthonormal basis of the
dual of the Cartan subalgebra of $C_n$ with
\begin{equation}\label{eq:croot}
\alpha'_0={1\over2}\delta-\epsilon'_1, \quad
\alpha'_1=\epsilon'_1-\epsilon'_2,\quad\dots,\quad
\alpha'_{n-1}=\epsilon'_{n-1}-\epsilon'_n,\quad
\alpha'_n=2\epsilon'_n.
\end{equation}
Comparing (\ref{eq:broot}) and (\ref{eq:croot}),
we have the relation
$\epsilon_i+\epsilon'_{n+1-i}={1\over2}\delta$.

The unnormalized character of ${\cal L}(\Lambda'_0)$
 is given by \cite{kacp}
\begin{equation*}
{\rm ch}\, {\cal L}(\Lambda'_0)=\left.\Theta\right/
\prod_{j=1}^\infty(1-e^{-j\delta})^n,\quad
\Theta=
\sum_{\gamma'\in M'}
e^{\Lambda'_0+\gamma' - {1\over2}
|\gamma'|^2\delta},
\quad M'=\bigoplus_{i=1}^n{\bf Z}\epsilon'_i.
\end{equation*}
The fundamental weight $\Lambda'_0$ in the $C_n$ picture
is related to $\Lambda_n$
in the $B_n$ picture as $\Lambda'_0=\Lambda_n+b_0\delta$ with
a certain constant $b_0$. Thus
$e^{-b_0\delta} {\rm ch}\, {\cal L}(\Lambda'_0)$ is
 the unnormalized character of
${\cal L}(\Lambda_n)$.
After the substitution of $\Lambda'_0=\overline\Lambda_n
+{1\over2}\Lambda_0+b_0\delta$
and $\epsilon'_i={1\over2}\delta-\epsilon_{n+1-i}$,
$e^{-b_0\delta}\Theta$ is expressed as
\begin{equation*}
e^{-b_0\delta}\Theta=e^{{1\over2}|\overline\Lambda_n|^2
\delta}
\sum_{\gamma\in M}
e^{{1\over2}\Lambda_0+\overline\Lambda_n+\gamma - {1\over2}
|\overline\Lambda_n+\gamma|^2\delta},
\quad M=\bigoplus_{i=1}^n{\bf Z}\epsilon_i,
\end{equation*}
where $\overline\Lambda_n={1\over2}(\epsilon_1+\cdots
+\epsilon_n)$.
Finally, by setting $e^{-\delta}=q$ and $e^{\Lambda_0}=1$,
we have
\begin{equation*}
{\rm ch}\, {\cal L}(\Lambda_n)(q,x)
={q^{-{1\over2}|\overline\Lambda_n|^2}\over
(q)_\infty^n}
\sum_{\gamma\in M} q^{{1\over2}|\overline\Lambda_n+
\gamma|^2}e^{\overline\Lambda_n+\gamma},\quad
e^{ \epsilon_i}=x_i,
\end{equation*}
which is the left hand side of (\ref{eq:atwocharacter}).

\end{document}